\documentclass[leqno]{article}
\usepackage{CJK,CJKnumb,CJKulem,times,dsfont,ifthen,mathrsfs,latexsym,amsfonts,color}
\usepackage{amsmath,amsthm,makeidx,fontenc,amssymb,bm,graphicx,psfrag,listings,curves,extarrows,enumitem}

\usepackage{geometry}
\usepackage{indentfirst} 
\usepackage{amssymb}
\makeatletter

\newcommand{\Rmnum}[1]{\expandafter\@slowromancap\romannumeral #1@}
\makeatother

\newtheorem{definition}{Definition}[section]
\newtheorem{theorem}{Theorem}[section]
\newtheorem{lemma}{Lemma}[section]
\newtheorem{corollary}{Corollary}[section]

\newtheorem{remark}{Remark}[section]

\newcommand{\al}{\alpha}
\newcommand{\ga}{\gamma}

\newcommand{\dl}{\Delta}
\newcommand{\e}{\varepsilon}
\newcommand{\om}{\Omega}
\newcommand{\iy}{\infty}

\newcommand{\la}{\lambda}
\newcommand{\vp}{\varphi}
\newcommand{\pa}{\partial}
\newcommand{\ti}{\tilde}

\newcommand{\ra}{\rightarrow}
\newcommand{\rh}{\rightharpoonup}
\newcommand{\hra}{\hookrightarrow}

\newcommand{\lab}{\label}
\newcommand{\f}{\frac}
\newcommand{\bt}{\begin{theorem}}
\newcommand{\et}{\end{theorem}}
\newcommand{\bl}{\begin{lemma}}
\newcommand{\el}{\end{lemma}}
\newcommand{\bd}{\begin{definition}}
\newcommand{\ed}{\end{definition}}
\newcommand{\bc}{\begin{corollary}}
\newcommand{\ec}{\end{corollary}}
\newcommand{\bp}{\begin{proof}}
\newcommand{\ep}{\end{proof}}
\newcommand{\bx}{\begin{example}}
\newcommand{\ex}{\end{example}}
\newcommand{\bi}{\begin{exercise}}
\newcommand{\ei}{\end{exercise}}
\newcommand{\br}{\begin{remark}}
\newcommand{\er}{\end{remark}}
\newcommand{\be}{\begin{equation}}
\newcommand{\ee}{\end{equation}}
\newcommand{\bal}{\begin{align}}
\newcommand{\bn}{\begin{enumerate}}
\newcommand{\en}{\end{enumerate}}
\newcommand{\ba}{\begin{align}}
\newcommand{\bg}{\begin{align*}}
\newcommand{\eg}{\end{align*}}
\newcommand{\bcs}{\begin{cases}}
\newcommand{\ecs}{\end{cases}}
\newcommand{\C}{{\mathbb C}}
\newcommand{\R}{{\mathbb R}}
\newcommand{\RN}{{\mathbb R^N}}
\newcommand{\bean}{\begin{eqnarray*}}
\newcommand{\eean}{\end{eqnarray*}}

\newcommand{\s}{\star}
\newcommand{\mc}{\mathcal}
\newcommand{\D}{{\mathcal D}}
\newcommand{\p}{{\mathcal P}}
\newcommand{\sbr}[1]{\left(#1\right)}
\newcommand{\mbr}[1]{\left[#1\right]}
\newcommand{\lbr}[1]{\left\{#1\right\}}

\setitemize{itemindent=38pt,leftmargin=0pt,itemsep=-0.4ex,listparindent=26pt,partopsep=0pt,parsep=0.5ex,topsep=-0.25ex}

\numberwithin{equation}{section}

\begin{document}
\theoremstyle{plain}

\title{\bf Normalized ground states for semilinear elliptic systems with critical and subcritical nonlinearities\thanks{This work is supported by NSFC(11801581,11025106, 11371212, 11271386);  E-mails: li-hw17@mails.tsinghua.edu.cn\quad\& \quad zou-wm@mail.tsinghua.edu.cn}}

\date{}
\author{
{\bf Houwang Li$^1$\;\&\;Wenming Zou$^2$}\\
\footnotesize \it 1. Department of Mathematical Sciences, Tsinghua University, Beijing 100084, China.\\
\footnotesize \it 2. Department of Mathematical Sciences, Tsinghua University, Beijing 100084, China.}

\maketitle
\begin{center}
\begin{minipage}{120mm}
\begin{center}{\bf Abstract }\end{center}
In the present paper, we study the normalized solutions with least energy to the following system:
$$\begin{cases}
-\Delta u+\la_1u=\mu_1 |u|^{p-2}u+\beta r_1|u|^{r_1-2}|v|^{r_2}u\quad&\hbox{in}~\RN,\\
-\Delta v+\la_2v=\mu_2 |v|^{q-2}v+\beta r_2|u|^{r_1}|v|^{r_2-2}v\quad&\hbox{in}~\RN,\\
\int_{\RN}u^2=a_1^2\quad\hbox{and}\quad\int_{\RN}v^2=a_2^2,
\end{cases}$$
where $p,r_1+r_2<2^*$ and $q\le2^*$. To this purpose, we study the geometry of the Pohozaev manifold
and the associated minimizition problem. Under some assumptions on $a_1,a_2$ and $\beta$,
we obtain the existence of the positive normalized ground state solution to the above system.
\vskip0.23in

{\bf Key  words:}   Semilinear elliptic system; Normalized ground states; Pohozaev manifold; Sobolev critical.

\vskip0.1in
{\bf 2010 Mathematics Subject Classification:} 35J50, 35J15, 35J60.

\vskip0.23in

\end{minipage}
\end{center}

\vskip0.26in
\newpage
\section{Introduction}

	We recall the following Schr\"odinger system:
\be\lab{201905-1}
\begin{cases}
-i\f{\pa}{\pa t}\Phi_1=\dl \Phi_1+\mu_1 |\Phi_1|^{p-2}\Phi_1+\beta r_1|\Phi_1|^{r_1-2}|\Phi_2|^{r_2}\Phi_1,\\
-i\f{\pa}{\pa t}\Phi_2=\dl \Phi_2+\mu_2|\Phi_2|^{q-2}\Phi_1+\beta r_2|\Phi_1|^{r_1}|\Phi_2|^{r_2-2}\Phi_2,\\
\Phi_j=\Phi_j(x,t)\in \C,\ (x,t)\in \RN\times \R,\ j=1,2,
\end{cases}
\ee
where $i$ is the imaginary unit, $\mu_1,\mu_2$ and $\beta$ are constants, which comes from various physical phenomena, such as mean-field modles for binary mixtures
of Bose-Einstein condensates, or binary gases of fermion atoms in degenerate quantum states (Bose-Fermi mixtures, Fermi-Fermi mixtures), 
see \cite{Ad.2007,BFKM.2015,EGBB.1997,Mal.2008} for more physical background.
Physically, system \eqref{201905-1} has the nature of conservation of mass, that is the following two norms
	$$\int_{\RN}|\Phi_1(t,x)|^2\mathrm{d}x\quad \text{and}\quad\int_{\RN}|\Phi_2(t,x)|^2\mathrm{d}x$$
are independent of $t\in\R$. Moreover, the $L^2$-norms $|\Phi_1(t,\cdot)|_2$ and $|\Phi_2(t,\cdot)|_2$ have important physical significance, for example,
 in Bose-Einstein condensates, $|\Phi_1(t,\cdot)|_2$ and $|\Phi_2(t,\cdot)|_2$
 represent the number of particles of each component; in nonlinear optics framwork, $|\Phi_1(t,\cdot)|_2$ and $|\Phi_2(t,\cdot)|_2$
 represent the power supply.
Therefore it is natural to consider the masses as preserved, and the solution of \eqref{201905-1} with prescribed mass is called normalized solution.

\vskip0.08in

	In order to study the solitary wave solution of \eqref{201905-1}, we set $\Phi_1(x,t)=e^{i\la_1 t}u(x)$ and $\Phi_2(x,t)=e^{i\la_2t}v(x)$.
Then the system \eqref{201905-1} is reduced to the general elliptic system:
\be\lab{201905-2}
	\begin{cases}
	-\Delta u+\la_1u=\mu_1 |u|^{p-2}u+\beta r_1|u|^{r_1-2}|v|^{r_2}u\quad&\hbox{in}\;\RN,\\
	-\Delta v+\la_2v=\mu_2 |v|^{q-2}v+\beta r_2|u|^{r_1}|v|^{r_2-2}v\quad&\hbox{in}\;\RN.
	\end{cases}
\ee
And the existence of normalized solutions to \eqref{201905-2} can be formulated as follows:
given $a_1,a_2>0$, we aim to find $(u,v)\in H^1(\RN)\times H^1(\RN)$ and $(\la_1,\la_2)\in\R^2$ such that
\be\lab{201905-3}
	\begin{cases}
	-\Delta u+\la_1u=\mu_1 |u|^{p-2}u+\beta r_1|u|^{r_1-2}|v|^{r_2}u\quad &\hbox{in}\;\RN,\\
	-\Delta v+\la_2v=\mu_2 |v|^{q-2}v+\beta r_2|u|^{r_1}|v|^{r_2-2}v\quad&\hbox{in}\;\RN,\\
	\int_{\RN}u^2=a_1^2\quad\text{and}\quad\int_{\RN}v^2=a_2^2.
	\end{cases}
\ee
Throughout the paper, we treat \eqref{201905-3} in cases $\mu_1,\mu_2,\beta>0$, which is the so-called self-focusing and attractive interaction, 
and we require also
\be\lab{2020-07zwm} 
	N\ge3,\quad  r_1,r_2>1, \quad  2<p<2^*,\quad   2<r_1+r_2<2^*,\quad  2<q\le 2^*,
\ee
where $2^*=\f{2N}{N-2}$ is the Sobolev critical exponent. These constants are prescribed while
the parameters $\la_1,\la_2$ are unknown. It is easy to see that a normalized solution of \eqref{201905-3} can be found as a critical 
point of the energy functional
\be\lab{201905-engfun}
	I(u,v)=\int_{\RN}\f{1}{2}(|\nabla u|^2+|\nabla v^2|)-\f{1}{p}\mu_1|u|^p-\f{1}{q}\mu_2|v|^p-\beta|u|^{r_1}|v|^{r_2}
\ee
under the constraint $S_{a_1}\times S_{a_2}$, where
	$$S_{a}=\left\{u\in H^1(\RN):\int_{\RN}u^2=a^2\right\},$$
and the parameters $\la_1,\la_2$ appear as Lagrangian multipliers. In this paper, we are particularly interested in the
normalized ground states
\bd
We say that $(u_0,v_0)$ is a normalized ground state of system \eqref{201905-3}, if it is a solution to \eqref{201905-3}
having minimal energy amoung all the normalized solutions:
	$$I(u_0,v_0)=\inf\lbr{I(u,v):(u,v)~\text{solves}~\eqref{201905-3}~\text{for some}~(\la_1,\la_2)\in\R^2}.$$
\ed

The search for normalized ground states of system \eqref{201905-3} is a challenging and interesting problem.
The presence of the $L^2$-constraint makes the methods developed to deal with unconstraint problems unavailable, and
new technical difficulties arise. One of the main difficulties is the lack of the compactness of the constraint Palais-Smale
sequences. Indeed it is hard to check that the weak limits of the constraint Palais-Smale sequences lie in the constraint
$S_{a_1}\times S_{a_2}$, since the embeddings $H^1(\RN)\hra L^2(\RN)$ and even $H^1_{rad}(\RN)\hra L^2(\RN)$ are not compact.
Moreover, the $L^2$-constraint induces  a new critical exponent, the $L^2$-critical exponent
	$$\bar p=2+\f{4}{N}.$$
This is the threshold exponent for the boundedness of the energy functional $I(u,v)$. If the problem is purely $L^2$-subcritical i.e.,
 $2<p,q,r_1+r_2<\bar p$, then $I(u,v)$ is bounded from below on $S_{a_1}\times S_{a_2}$. In this case, T. Gou and L. Jeanjean 
\cite{GouJeanjean=2016} obtained the compactness of the minimizing sequence of $I(u,v)$ constrianed on $S_{a_1}\times S_{a_2}$,
and proved the existence of a normalized ground state as a global minimizer. However, if one of $p,q,r_1+r_2$ is greater
than $\bar p$, i.e., $L^2$-supercritical, then $I(u,v)$ is unbounded from below and from above on $S_{a_1}\times S_{a_2}$. In the cases
 $2<p,q<\bar p<r_1+r_2<2^*$ and $2<r_1+r_2<\bar p<p,q<2^*$,  T. Gou
and L. Jeanjean \cite{GouJeanjean=2018} proved the existence of a normalized ground state; in the cases $\bar p<p,q,r_1+r_2<2^*$,
using the Pohozaev manifold and mountain pass lemma, T. Bartsch, L. Jeanjean and N. Soave  (See  \cite{BJS=JMPA=2016,BartschJeanjean=2018}) proved the existence of a normalized ground state for
large $\beta$, below we will  give more detailed comparisons on these results.  For  other  conclusions about the existence and mulplicity of the normalized
solutions for Schr\"odinger equations on the whole space, we refer to \cite{BartschSoave=2017,GouJeanjean=2018,BJS=JMPA=2016,
BartschJeanjean=2018,BS=CVPDE=2019,Soave=arXiv=2018,Soave=arXiv=2019,Ikoma=arXiv=2019}.

\vskip0.12in

	We note that in \cite{Soave=arXiv=2019},
N. Soave considered the following  nonlinear Schr\"{o}dinger equation with combined power nonlinearities:
\be\label{zwm78}
	-\Delta u = \lambda u + \mu |u|^{q-2} u + |u|^{2^\ast-2}u \quad \text{in} \; \R^N,\;  N\geq 3, 
\ee
with  prescribed mass
 $$ \int_{\R^N} |u|^2 = a^2,$$
in the Sobolev critical case. For a $L^2$-subcritical, $L^2$-critical, and $L^2$-supercritical perturbation $\mu |u|^{q-2}u$, the author  proved several existence/non-existence and stability/instability results.  He obtained a constraint Palais-Smale sequence with an additional property by studying the geometry of the corrsponding
Pohozaev manifold, and he proved the compactness of this special constraint Palais-Smale sequence under some energy level.
We are motivated by  \cite{Soave=arXiv=2019}  to study  the system \eqref{201905-3}. However, we deal with a system, which is different from the   scalar equation: the appearence of the coupled item makes the geometry of the Pohozaev manifold more complicated; the
compactness of constraint Palais-Smale sequence is hard  to  get.

\vskip0.1in

	For simplicity, let $r=r_1+r_2$ and
\be\lab{201905-gam}
	\ga_p=\f{N(p-2)}{2p}\left\{
		\begin{aligned}
			&<\f{2}{p}, &\text{if}~&2<p<\bar p,\\	 &=\f{2}{p}, &\text{if}~&p=\bar p,\\	&>\f{2}{p}, &\text{if}~&\bar p<p<2^*,
		\end{aligned}\right.
		\quad\text{and}\quad\ga_{2^*}=1.
\ee
As in \cite{Soave=arXiv=2018,Soave=arXiv=2019}, the following Pohozaev manifold will play a special role in the proof:
\be\lab{201905-phoman}
\mathcal P_{a_1,a_2}=\left\{(u,v)\in S_{a_1}\times S_{a_2}:P(u,v)=0\right\},\ee
where
\be\lab{201905-phoequ}
P(u,v)=\int_{\RN}|\nabla u|^2+|\nabla v|^2-\ga_p\mu_1|u|^p-\ga_q\mu_2|v|^q-r\ga_r\beta|u|^{r_1}|v|^{r_2}.\ee
As a consequence of the Pohozaev identity, any solution of \eqref{201905-3} belongs to $\mathcal P_{a_1,a_2}$. So
if $(u,v)\in\mc P_{a_1,a_2}$ is a minimizer of the constraint minimization
\be\lab{201905-min}
m(a_1,a_2)=\inf_{(u,v)\in\mc P_{a_1,a_2}}I(u,v),\ee
and $(u,v)$ solves system \eqref{201905-2} for some $\la_1,\la_2$, then $(u,v)$ is a normalized ground state of \eqref{201905-3}.
To study the minimization problem \eqref{201905-min}, we introduce a dilition operation preserving the $L^2$-norm:
for $u\in S_a$ and $s\in\R$,
	$$s\s u(x):=e^{\f{Ns}{2}} u(e^sx)\quad \text{for a.e.}~x\in\RN.$$
Then $s\s u\in S_a$. Define $s\s(u,v)=(s\s u,s\s v)$ and the fiber maps
\be\lab{201905-fibmap}
	\begin{aligned}
		\Phi_{(u,v)}(s):&=I(s\s (u,v))\\&=\int_{\RN}\f{e^{2s}}{2}(|\nabla u|^2+|\nabla v|^2)-\f{e^{p\gamma_ps}}{p}\mu_1|u|^p-\f{e^{q\gamma_qs}}{q}
			\mu_2|v|^q-e^{r\gamma_rs}\beta|u|^{r_1}|v|^{r_2}.
	\end{aligned}
\ee
By direct computation, we have $\Phi_{(u,v)}'(s)=P(s\s(u,v))$ and 
	$$\mc P_{a_1,a_2}=\left\{ (u,v)\in S_{a_1}\times S_{a_2}:\Phi_{(u,v)}'(0)=0 \right\}.$$
In this direction, we decompose $\p_{a_1,a_2}$ into disjoint unions $\p_{a_1,a_2}=\p^+_{a_1,a_2}\cup\p^0_{a_1,a_2}\cup\p^-_{a_1,a_2}$,
where
\begin{equation*}\begin{aligned}
	\mc P_{a_1,a_2}^+&:=\left\{ (u,v)\in S_{a_1}\times S_{a_2}:\Phi_{(u,v)}''(0)>0\right\},\\
	\mc P_{a_1,a_2}^0&:=\left\{ (u,v)\in S_{a_1}\times S_{a_2}:\Phi_{(u,v)}''(0)=0\right\},\\
	\mc P_{a_1,a_2}^-&:=\left\{ (u,v)\in S_{a_1}\times S_{a_2}:\Phi_{(u,v)}''(0)<0\right\}.
\end{aligned}\end{equation*}
We see that the monotonicity and convexity of $\Phi_{(u,v)}(s)$ will strongly affect the structure of $\mathcal P$ and
hence have a strong impact on the minimization problem \eqref{201905-min}.

\vskip0.08in

	Now, we state our main results. As we have stated, throughout this paper, we require $\mu_1,\mu_2,\beta,a_1,a_2>0$ and $r_1,r_2>1$.
For the convenience of description, let
\be
	T(a_1,a_2)=\left\{\begin{aligned}
 		&a_1^{r_1(1-\ga_r)}a_2^{r_2(1-\ga_r)}\beta(\mu_2a_2^{q(1-\ga_q)})^{\f{2-r\ga_r}{q\ga_q-2}}+\mu_1a_1^{p(1-\ga_p)}
 		(\mu_2a_2^{q(1-\ga_q)})^{\f{2-p\ga_p}{q\ga_q-2}}, 								&\text{if}~r&<\bar p,\\
		&\min\lbr{a_1^{r_1(1-\ga_r)}a_2^{r_2(1-\ga_r)}\beta,(\mu_1a_1^{p(1-\ga_p)})^{\f{1}{2-p\ga_p}}(\mu_2a_2^{q(1-\ga_q)})^
		{\f{1}{q\ga_q-2}}}, 															&\text{if}~r&=\bar p,\\
		&a_1^{r_1(1-\ga_r)}a_2^{r_2(1-\ga_r)}\beta(\mu_1a_1^{p(1-\ga_p)})^{\f{r\ga_r-2}{2-p\ga_p}}+\mu_2a_2^{q(1-\ga_q)}
		(\mu_1a_1^{p(1-\ga_p)})^{\f{q\ga_q-2}{2-p\ga_p}}, 								&\text{if}~r&>\bar p.
	\end{aligned}\right.
\ee
	Then we have a result concerning a mixed situation.
\bt\lab{thm2}
	Suppose $3\le N\le4$, $2<p<\bar p<q\le2^*$, $r<2^*,r_2<2$, then there exists a constant $\al_0=\al_0(p,q,r,N)>0$ such that
	if $T(a_1,a_2)<\al_0$, then \eqref{201905-3} has a positive normalized ground state.
\et
\br
	The assumption $r_2<2$ is used to control the energy level, and the assumption $T(a_1,a_2)<\al_0$ is applied to ensure that
	the Pohozaev manifold has a good geometry. We note that for fixed $\mu_1,\mu_2,\beta>0$, the constant  $T(a_1,a_2)<\al_0$ holds as long
	as $a_1a_2$ small enough.
\er

 	We also obtain a result about the normalized ground state for purely $L^2$-supercritical case.
\bt\lab{thm4}
	Suppose $3\le N\le4$, $\bar p<p,\ q,\ r<2^*$, then
	\begin{itemize}[fullwidth,itemindent=0em]
	\item[(1)]there exists a $\beta_0>0$ such that \eqref{201905-3} has a positive normalized ground state for any $\beta>\beta_0$;
	\item[(2)]if further $r_1,r_2<2$, then \eqref{201905-3} has a positive normalized ground state for any $\beta>0$.
	\end{itemize}
\et
\br
	The first conclusion of Theorem \ref{thm4} is similar to the results in \cite{BartschJeanjean=2018,BJS=JMPA=2016},
	but the second one is a new result. Particularly, let $N=4,\ p=q=2r_1=2r_2\in(\bar p,2^*)$, then according to the second conlusion 
	of Theorem \ref{thm4}, the following system 
	\begin{equation*}
		\begin{cases}
		-\Delta u+\la_1u=\mu_1 |u|^{p-2}u+\beta r_1|u|^{\f{p}{2}-2}|v|^{\f{p}{2}}u\quad &\hbox{in}~\R^4,\\
		-\Delta v+\la_2v=\mu_2 |v|^{p-2}v+\beta r_2|u|^{\f{p}{2}}|v|^{\f{p}{2}-2}v\quad &\hbox{in}~\R^4,\\
		\int_{\R^4}u^2=a_1^2,\quad\int_{\R^4}v^2=a_2^2,
		\end{cases}
	\end{equation*}
	has a positive normalized ground state for any $a_1,a_2,\beta>0$.
\er

The paper is organized as follows. In Section 2 we collect some preliminary results which will be used from time to time in the paper.
Theorems \ref{thm2}, \ref{thm4} are proved
in Sections 4, 5 respectively. In Appendix B, we give a proof of a regularity result. Throughout the paper we use the notation $|u|_{p}$
to denote the $L^p(\RN)$ norm, and we simply write $H^1=H^1(\RN),H=H^1(\RN)\times H^1(\RN)$. Similarly, $H_r^1$ denotes the subspace
of funtions in $H^1$  which are radial symmetric with respect to 0, and $H_r=H^1_r\times H^1_r,S_{a,r}=S_a\cap H^1_r$. The
symbol $||\cdot||$ denotes the norm in $H^1$ or $H$. Denoting by $u^*$ the symmetric decreasing rearrangement of $u\in H^1$, we
recall that (see \cite{LossLieb}) for $p,q>1$
	$$|\nabla u^*|_2\le|\nabla u|_2,\quad |u^*|_p=|u|_p \quad and\quad\int_{\RN}|u^*|^p|v^*|^q\ge\int_{\RN}|u|^p|v|^q.$$
The capital letters $C_1,C_2,\cdots$ denote positive constants which may depend on $N, $ $p,$ $q,$ $r_1,$ $r_2$, whose precise values can change
from line to line.

\vskip0.23in
\section{Preliminaries}

	In this section, we summarize several results which will be used in the rest disscussion.

\vskip0.12in
	For $N\geq3,2<p\le2^*$, the Gagliardo-Nirenberg inequality is
\be\lab{201905-GNinq}
|u|_p\leq \mc C_{N,p}|\nabla u|_2^{\ga_p}|u|_2^{1-\ga_p},\quad \forall u\in H^1,
\ee
where $\ga_p$ is defined by \eqref{201905-gam}. For a special case of \eqref{201905-GNinq}, if $p=2^*$,
then denoting $\mc S=\mc C_{N,2^*}^{-2}$, we have the Sobolev inequality
	$$\mc S|u|_{2^*}^2\leq |\nabla u|_2^2,\quad \forall u\in D^{1,2}(\RN),$$
where $D^{1,2}(\RN)$ is the completion of $C_c^{\iy}(\RN)$ with respect to the norm $||u||_{D^{1,2}}:=|\nabla u|_2$. We
observe that the functional $I(u,v)$ defined in \eqref{201905-engfun} is well defined and is of class $C^1$. Throughout
this paper, we denote
\be\lab{201905-conD}
	\bcs
	\D_1=(\f{\max\{r_1,r_2\}}{r})^{\f{r\ga_r}{2}}\mc C_{N,r}^ra_1^{r_1(1-\ga_r)}a_2^{r_2(1-\ga_r)},\\
	\D_2=\f{1}{p}\mu_1\mc C_{N,p}^{p}a_1^{p(1-\ga_p)},\\
	\D_3=\f{1}{q}\mu_2\mc C_{N,q}^{q}a_2^{q(1-\ga_q)}.
\ecs\ee
Then we have
\be\lab{201905-coupled}
	\begin{aligned}
	\quad\ \int|u|^{r_1}|v|^{r_2}&\leq(\int|u_1|^r)^{\f{r_1}{r}}(\int|u_2|^r)^{\f{r_2}{r}}\\
	&\leq \mc C_{N,r}a_1^{r_1(1-\ga_r)}a_2^{r_2(1-\ga_r)}(\int|\nabla u|^2)^{\frac{r_1\ga_r}{2}}
		(\int|\nabla v|^2)^{\f{r_2\ga_r}{2}}\\
	&\leq \mc C_{N,r}a_1^{r_1(1-\ga_r)}a_2^{r_2(1-\ga_r)}(\f{r_1}{r}\int|\nabla u|^2+\f{r_2}{r}\int|\nabla v|^2)^
		{\f{r\ga_r}{2}}\\
	&\leq \D_1(\int|\nabla u|^2+|\nabla v|^2)^{\f{r\ga_r}{2}},
	\end{aligned}
\ee
\be\lab{201905-nonlinearity}
	\f{1}{p}\int\mu_1|u|^p\leq \D_2|\nabla u|_2^{p\ga_p}\quad \text{and}\quad \f{1}{q}\int\mu_2|u|^q\leq\D_3|\nabla v|_2^{q\ga_q}.
\ee
Substituting \eqref{201905-coupled}-\eqref{201905-nonlinearity} into \eqref{201905-engfun}, we obtain
\be\lab{201911-11}
	\begin{aligned}
	I(u,v)&\geq\f{1}{2}(\int|\nabla u|^2+|\nabla v|^2)-\D_1\beta(\int|\nabla u|^2+|\nabla v|^2)^{\f{r\ga_r}{2}}
			-\D_2|\nabla u|_2^{p\ga_p}-\D_3|\nabla v|_2^{q\ga_q} \\
		&\geq h((\int|\nabla u|^2+|\nabla v|^2)^{\f{1}{2}}),
	\end{aligned}
\ee
where $h(t):(0,+\iy)\ra\R$ defined by
\be\lab{201911-h}
	h(t)=\f{1}{2}t^2-\D_1\beta t^{r\ga_r}-\D_2t^{p\ga_p}-\D_3t^{q\ga_q}.
\ee

\vskip0.12in

	We now focus on the Sobolev subcritical and critical nonlinear Schr\"odinger equations with prescribed $L^2$-norm. For
fixed $a>0,\mu>0,2<p\le2^*$, we search for $(u,\la)\in H^1\times\R$ solving
\be\lab{201905-nlse}
	\bcs
	-\dl u+\la u=\mu |u|^{p-2}u \quad \text{in}\ \RN,\\
	\int_{\RN}u^2 =a^2,\quad u\in H^1.
	\ecs
\ee
Solutions of \eqref{201905-nlse} can be found as the critical points of $E_{p,\mu}:H^1\to\R$
	$$E_{p,\mu}(u)=\int_{\RN}\f{1}{2}|\nabla u|^2-\f{1}{p}\mu |u|^p,$$
constrained on $S_a$, and the parameter $\la$ appears as Lagrangian multiplier. Back to the problem \eqref{201905-nlse}, testing the
 equation with $u$, we get
 	$$\int_{\RN}|\nabla u|^2+\la \int_{\RN} u^2=\mu \int_{\RN} |u|^p,$$
and then combined with the Pohozaev identity
	$$\f{N-2}{2}\int_{\RN}|\nabla u|^2+\f{N}{2}\la \int_{\RN} u^2=\f{N}{p}\mu \int_{\RN} |u|^p,$$
it follows that
\be\lab{202006-01}
	\la \int_{\RN} u^2=(1-\ga_p)\mu \int_{\RN} |u|^p.
\ee
Hence,  if $p<2^*$, then we must have  $\la>0$. Therefore,  by scaling, the equation \eqref{201905-nlse} is equivalent to
\be\lab{201905-w}
	-\dl w+w=|w|^{p-2}w\quad \text{in}\ \RN,\quad w\in H^1.
\ee
While if $p=2^*$, we have $\la=0$, then by scaling, the equation \eqref{201905-nlse} is equivalent to
\be\lab{201905-wc}
	-\dl w=|w|^{2^*-2}w\quad \text{in}\ \RN,\quad w\in D^{1,2}(\RN).
\ee
Since the  positive solutions of \eqref{201905-w}\eqref{201905-wc} are studied clearly, the existence of normalized solutions 
of \eqref{201905-nlse} can be obtained by scaling. However, there are still some special properties that need to be clarified. 
To be precise, we introduce the Pohozaev manifold for single equations
\be\lab{201905-phomans}
	\mc T_{a,p,\mu}:=\left\{u\in S_a:\int_{\RN}|\nabla u|^2-\ga_p\mu |u|^p=0\right\},
\ee
and the constraint minimizition problem
\be\lab{201905-mins}
	m_p^\mu(a)=\inf_{u\in\mc T_{a,p,\mu}}E_{p,\mu}(u).
\ee
It is easy to see that
	$$m(a_1,0)=m_p^{\mu_1}(a_1)\quad \text{and}\quad m(0,a_2)=m_q^{\mu_2}(a_2).$$
We have the following lemmas.

\bl\lab{201905-nlses}
	Suppose $N\ge3,\mu, a>0$ and $2<p<2^*, p\neq \bar p$, then up to a translation, \eqref{201905-nlse} has a unique positive
	solution $u_{p,\mu}\in \mc T_{a,p,\mu}$ with $\la>0$. Moreover,
	\begin{itemize}[fullwidth,itemindent=0em]
		\item[(1)]	if $p<\bar p$, then
		\be\lab{201905-subc}
			m_p^\mu(a)=\inf_{u\in S_a}E_{p,\mu}(u)=E_{p,\mu}(u_{p,\mu})<0;
		\ee
		\item[(2)]	if $p>\bar p$, then
		\be\lab{201905-superc}
			m_p^\mu(a)=\inf_{u\in S_a}\max_{t\in\R}E_{p,\mu}(t\s u)=\max_{t\in\R}E_{p,\mu}(t\s u_{p,\mu})=E_{p,\mu}(u_{p,\mu})>0;
		\ee
	\end{itemize}
	Finally,  for  both cases $m_p^\mu(a)$ is strictly decreasing with respect to $a>0$.
\el
\bp
By \cite{Kwong=1989,GNN=1979}, up to a translation, $w_{p,\mu}$ is the unique positive solution of \eqref{201905-w}, which is
radial symmetric and decreasing with respect to 0. Since $p\ne\bar p$, by scaling we obtain the unique solution of \eqref{201905-nlse}
	$$u_{p,\mu}=(\f{\la}{\mu})^{\f{1}{p-2}}w_{p,\mu}(\la^{\f{1}{2}}x)	\quad \text{with} \quad
			\la=(\f{a^2}{|w_{p,\mu}|_2^2}\mu^{\f{2}{p-2}})^{\f{p-2}{2-p\ga_p}}.$$
Using the Pohozaev identity, it is easy to check that $u_{p,\mu}\in\mc T_{a,p,\mu}$. On the other hand, the equality
of \eqref{201905-GNinq} is achieved by $w_{p,\mu}$, that is
\be\lab{202006-02}
	|w_{p,\mu}|_p=\mc C_{p,N}|\nabla w_{p,\mu}|_2^{\ga_p}|w_{p,\mu}|_2^{1-\ga_p},
\ee
and taking a similar process as the one to get equation \eqref{202006-01}, we obtain 
\be\lab{202006-03}
	|\nabla w_{p,\mu}|_2^2=\ga_p|w_{p,\mu}|_p^p.
\ee
Then combining \eqref{202006-02} and \eqref{202006-03}, there is $|\nabla w_{p,\mu}|_2^{2-p\ga_p}=\ga_pC_{p,N}|w_{p,\mu}|_2^{p-p\ga_p}$. So
\be\lab{202006-04}
	\begin{aligned}
	E_{p,\mu}(u_{p,\mu})&=(\frac{1}{2}-\f{1}{p\ga_p})\int_{\RN}|\nabla u_{p,\mu}|^2\\
		&=(\frac{1}{2}-\f{1}{p\ga_p})\mu^{-\f{2}{p-2}}\la^{\f{p-p\ga_p}{p-2}}|\nabla w_{p,N}|_2^2\\
		&=(\frac{1}{2}-\f{1}{p\ga_p})(\mu a^{p-p\ga_p})^{\f{2}{2-p\ga_p}}\f{|\nabla w_{p,N}|_2^2}{|w_{p,N}|_2^{\f{2(p-p\ga_p)}{2-p\ga_p}}}\\
		&=(\frac{1}{2}-\f{1}{p\ga_p})(\ga_p\mc C_{N,p}\mu a^{p-p\ga_p})^{\f{2}{2-p\ga_p}},
	\end{aligned}
\ee
which is negative if $p<\bar p$ and is positive if $p>\bar p$. To prove futher properties, let
\be
	\begin{aligned}
	\Phi_u(s)&:=E_{p,\mu}(s\s u)\\
		&=\f{e^{2s}}{2}\int_{\RN}|\nabla u|^2-\f{e^{p\ga_ps}}{p}\ga_p\int_{\RN}\mu|u|^p.
	\end{aligned}
\ee
For the case $p<\bar p$, we take $u\in S_a$, then a direct computation tells us that there exists a unique global minimizer $s_u$ 
for $\Phi_u(s)$ and $s_u\s u\in\mc T_{a,p,\mu}$. So
	$$E_{p,\mu}(u)\ge E_{p,\mu}(s_u\s u)\ge m_p^\mu(a)\ge\inf_{u\in S_a}E_{p,\mu}(u),$$
which implies $m_p^\mu(a)=\inf_{u\in S_a}E_{p,\mu}(u)<0$. Taking a minimizing sequence $u_n$ for $\inf_{u\in S_a}E_{p,\mu}(u)$,
we assume $u_n\ge0,~u_n\in H^1_r$ by insteading $u_n$ of $|u_n|^*$. The coerciveness of $E_{p,\mu}|_{S_a}$ means
that $u_n$ is bounded. Then up to a subsequence $u_n\rh u_0$ in $H^1(\RN)$, $u_n\ra u_0$ in $L^p(\RN)$ and $u_n\ra u_0$ a.e. in $\RN$.
So $u_0\ge0$. We will prove that $u_0$ is a nontivial minimizer of $m_p^\mu(a)$. If $u_0=0$, then we have
	$$m_p^\mu(a)=\lim_{n\ra\iy}E_{p,\mu}(u_n)=\lim_{n\ra\iy}\frac{1}{2}\int_{\RN}|\nabla u_n|^2\ge0,$$
in contradiction with $m_p^\mu(a)<0$. Hence, $0<|u_0|_2\le a$. Suppose $|u_0|_2\ne a$, then
	$$m_p^\mu(a)=E_{p,\mu}(u_n)+o(1)\ge E_{p,\mu}(u_0)\ge m_p^\mu(|u_0|_2).$$
On the other hand, following the proof of \cite[Lemma 3.1]{GouJeanjean=2016}, we have
	$$m_p^\mu(a)\le m_p^\mu(|u_0|_2)+m_p^\mu(a-|u_0|_2)<m_p^\mu(|u_0|_2),$$
which is a contradiction. Thus $|u_0|_2=a$, $E_{p,\mu}(u_0)=m_p^\mu(a)$ and $u_0$ is a positive solution of
\eqref{201905-nlse}. Then the uniqueness result implies $u_0=u_{p,\mu}$ and hence $m_p^\mu(a)=E_{p,\mu}(u_{p,\mu})$. Recall \eqref{202006-04},
we know that $m_p^\mu(a)$ is strictly decreasing with respect to $a>0$.

	Suppose now $p>\bar p$, then from \cite[Lemma 2.10]{Jeanjean=1997}, there is $m_p^\mu(a)=\max_{t\in\R}E_{p,\mu}(t\s u_{p,\mu})=E_{p,\mu}(u_{p,\mu})$, 
and we can prove $m_p^\mu(a)=\inf_{u\in S_a}\max_{t\in\R}E_{p,\mu}(t\s u)$ similarly as\cite[Proposition 2.2]{Soave=arXiv=2019}.
\ep

	For the Sobolev critical case $p=2^*$, there is also a clear characterization about the positive solutions of \eqref{201905-nlse} 
and the minimizition problem \eqref{201905-mins}.
\bl\lab{201905-nlsec}
	Suppose $N\ge3,\mu,a>0$ and $p=2^*$, then
	\be\lab{201905-sobcri}
		m_{2^*}^\mu(a)=\inf_{u\in S_a}\max_{t\in\R}E_{2^*,\mu}(t\s u)=\f{1}{N}\mu^{-\f{N-2}{2}}\mc S^{\f{N}{2}}>0.
	\ee
	Moreover,
	\begin{itemize}[fullwidth,itemindent=0em]
	\item[(1)]	if $N=3,4$, then \eqref{201905-nlse} has no posotive solution for any $\la\in\R$, and in particular $m_{2^*}^\mu(a)$ is not achieved;
	\item[(2)]	if $N\ge5$, then up to a translation, \eqref{201905-nlse} has a unique positive solution $u_{2^*,\mu}\in \mc T_{a,2^*,\mu}$ with $\la=0$, and
		$$m_{2^*}^\mu(a)=E_{2^*,\mu}(u_{2^*,\mu}).$$
	\end{itemize}
\el
\bp
For detailed proof, we refer to \cite[Propesition 2.2]{Soave=arXiv=2019}.
\ep

\vskip0.12in

	In the following paper, we need the following result.
\bl\lab{201905-liouville}
	Suppose $(u,v)\in H$ is a nonnegative solution of \eqref{201905-2} with $2<p,q,r\leq 2^*$, then
	\begin{itemize}[fullwidth,itemindent=0em]
	\item[(1)]	if $N=3,4$, then $u>0$ implies $\la_1>0$; $v>0$ implies $\la_2>0$;
	\item[(2)]	if $N\ge5$, then $u>0$ implies $\la_1\ge0$; $v>0$ implies $\la_2\ge0$.
	\end{itemize}
\el
\bp
From Corollary \ref{201905-Ac}, we know that $(u,v)$ is a smooth solution. Suppose $u>0$ but $\la_1<0$, then
	$$-\dl u=|\la_1|u+\mu_1 u^{p-1}+\beta r_1u^{r_1-1}v^{r_2}\ge \min\{|\la_1|,\mu_1\}u^\sigma \quad \text{in}\ \RN,$$
for any $1<\sigma<p-1$. Using a Liouville type theorem \cite[Theorem 8.4]{QuittnerSouplet}, we deduce $u=0$, which is impossible since $u>0$.
So $\la_1\ge0$. Morevoer, if $N=3,4$ and $\la_1=0$, i.e.,
	$$-\dl u=\mu_1 u^{p-1}+\beta r_1u^{r_1-1}v^{r_2}\ge0 \quad \text{in}\ \RN,$$
then \cite[Lemma A.2]{Ikoma=AdvNonStu=2014} implies that $u=0$, which is also a contradiction. So $\la_1>0$ when $N=3,4$.
\ep

	Finally we recall a Brezis-Lieb type lemma.
\bl\lab{201905-BL}
	Suppose $(u_n,v_n)\subset H$ is a bounded sequence, $(u_n,v_n)\to(u,v)$ a.e. in $\RN$ and $2\le r\le 2^*,r_1,r_2>1$, then
	$$\lim_{n\ra\iy}\int_{\RN}|u_n|^{r_1}|v_n|^{r_2}-|u|^{r_1}|v|^{r_2}-|u_n-u|^{r_1}|v_n-v|^{r_2}=0.$$
\el

\vskip0.23in
\section{The mixed exponent case}

	In this section, we assume $2<p<\bar p<q\le2^*, 2<r<2^*,r_1,r_2>1$. Recall the definition of $h(t)$ in \eqref{201911-h}, we have
\bl\lab{201905-4h}
	There exists a constant $\al_1>0$ such that if $T(a_1,a_2)<\al_1$, then the function $h(t)$ has exactly two critical points,
	one is a local minimum at negitive level, the other one is a global maximum at positive level. Futher, there exists $0<R_0<R_1$
	such that $h(R_0)=h(R_1)=0$, and $h(t)>0$ iff $t\in(R_0,R_1)$.
\el
\br\lab{tem1}
	Let $l(t)=at^2-bt^{r\ga_r}-ct^{p\ga_p}-dt^{q\ga_q}$ with $2<p<\bar p<q\le2^*$, $2<r<2^*$ and $a,b,c,d>0$. Then $l(t)$ has 
at most two critical points in $(0,\iy)$.
\er\lab{Remark1}
	The proof of Lemma \ref{201905-4h} and Remark \ref{tem1} is elementary but technique, we postpone it to the Appendix .

\bl\lab{201905-4pho}
	There exists a constant $\al_2>0$ such that if  $T(a_1,a_2)<\al_2$, then $\mc P^0_{a_1,a_2}=\emptyset$,
	and $\p_{a_1,a_2}$ is a $C^1$ submanifold in $H$ with codimension 3.
\el
\bp
	We first prove that $\p^0_{a_1,a_2}=\emptyset$ implies that $\p_{a_1,a_2}$ is a $C^1$ submanifold in $H$ with codimension 3. As
we can see, $\p_{a_1,a_2}$ is defined by $P(u,v)=0,G(u)=0,F(v)=0$, where
	$$G(u)=a_1^2-\int_{\RN} u^2,\quad F(v)=a_2^2-\int_{\RN} v^2.$$
It is sufficient to prove
	$$d(P,G,F):H\ra\R^3\quad \text{is a surjective}.$$
Suppose it is not true, by the independence of $dG(u)$ and $dF(v)$, there must be that $dP(u,v)$ is a linear combination of $dG(u)$
and $dF(v)$, i.e., there exist  $\nu_1,\nu_2\in\R$ such that $(u,v)$ is a weak solution of
\be\lab{201905-8}
	\bcs
	-\Delta u+\nu_1u=\f{p\ga_p}{2}\mu_1 |u|^{p-2}u+\f{r\ga_r}{2}\beta r_1|u|^{r_1-2}|v|^{r_2}u\quad \text{in}\  \RN,\\
	-\Delta v+\nu_2v=\f{q\ga_q}{2}\mu_2	|v|^{q-2}v+\f{r\ga_r}{2}\beta r_2|u|^{r_1}|v|^{r_2-2}v\quad \text{in}\  \RN,\\
	|u|_2=a_1,\quad |v|_2=a_2.
	\ecs
\ee
Testing system \eqref{201905-8} with $(u,v)$ and combining with the Pohozaev identity, we can conclude that
	$$2\int_{\RN}|\nabla u|^2+|\nabla v|^2=p\ga_p^2\int_{\RN}\mu_1|u|^p+q\ga_q^2\int_{\RN}\mu_2|v|^q+
		(r\ga_r)^2\int_{\RN}\beta |u|^{r_1}|v|^{r_2},$$
which implies that $(u,v)\in\p^0_{a_1,a_2}$, a contradiction.

Now we prove that there exists a constant $\al_2>0$ such that $\mc P^0_{a_1,a_2}=\emptyset$ as long as $T(a_1,a_2)<\al_2$.
Suppose there is a $(u,v)\in\mc P^0_{a_1,a_2}$. Let $\rho=(|u|_2^2+|v|_2^2)^{\f{1}{2}}$ and
	$$\begin{aligned}
	W(t):&=t\Phi_{(u,v)}'(0)-\Phi_{(u,v)}''(0)\\
		&=(t-2)\int_{\RN}|\nabla u|^2+|\nabla v|^2-(t-p\ga_p)\ga_p\int_{\RN}\mu_1 |u|^p\\
		&\quad-(t-q\ga_q)\ga_q\int_{\RN}\mu_2|v|^q-(t-r\ga_r)r\ga_r\int_{\RN}\beta|u|^{r_1}|v|^{r_2}\\
		&=0.
	\end{aligned}$$

\vskip0.12in
	We divide the proof into four different situations.
\begin{itemize}[fullwidth,itemindent=0em]
	\item[{\textit Case 1: $p\le r<\bar p$.}]	There is $p\ga_p\le r\ga_r<2<q\ga_q$. On the one hand, $W(r\ga_r)=0$ implies
		$$(2-r\ga_r)\rho^2\le(q\ga_q-r\ga_r)\int_{\RN}\mu_2|v|^q\le(q\ga_q-r\ga_r)q\D_3\rho^{q\ga_q}.$$
	It follows that $\rho\ge (\f{q\ga_q-r\ga_r}{2-r\ga_r}\f{1}{q\D_3})^{\f{1}{q\ga_q-2}}$. On the other hand, by $W(q\ga_q)=0$, we obtain
		\begin{align*}
			(q\ga_q-2)&=(q\ga_q-p\ga_p)\ga_p\rho^{-2}\int_{\RN}\mu_1 |u|^p+(q\ga_q-r\ga_r)r\ga_r\rho^{-2}\int_{\RN}\beta|u|^{r_1}|v|^{r_2}\\
				&\le(q\ga_q-p\ga_p)\ga_pp\D_2\rho^{p\ga_p-2}+(q\ga_q-r\ga_r)r\ga_r\D_1\beta\rho^{r\ga_r-2}\\
				&\le C(p,q,r)(\D_2\D_3^{\f{2-p\ga_p}{q\ga_q-2}}+\D_1\beta\D_3^{\f{2-r\ga_r}{q\ga_q-2}}).
		\end{align*}
	Therefore by the definitions of $\D_1,\D_2,\D_3$, we can choose an $\al_2>0$ such that
		$$\D_2\D_3^{\f{2-p\ga_p}{q\ga_q-2}}+\D_1\beta\D_3^{\f{2-r\ga_r}{q\ga_q-2}}<C(p,q,r)^{-1}(q\ga_q-2)$$
	as long as  $T(a_1,a_2)<\al_2$, then we get a contradiction. That is, $\mc P^0_{a_1,a_2}=\emptyset$ provided that $T(a_1,a_2)<\al_2$.

	\vskip0.08in
	\item[{\textit Case 2: $r<p<\bar p$.}]	If we exchange the roles played by $\D_2t^{p\ga_p}$ and $\D_1t^{r\ga_r}$ in {\it Case 1}, 
	then we can get a constant $\al_2$ with the required properties.

	\vskip0.08in
	\item[{\textit Case 3: $r=\bar p$.}]	We first suppose $\al_2<\frac{1}{4}$, so that $\f{1}{2}-\D_1\beta\in(\f{1}{4},\f{1}{2})$ when $\D_1\beta<\al_2$.
	Then completely analogous as {\it Case 1}, by combining $W(q\ga_q)=0$ and $W(p\ga_p)=0$, we can obtain the constant $\al_2$ with the required properties.

	\vskip0.08in
	\item[{\textit Case 4: $r>\bar p$.}]
	If $r\le q$, then there is $p\ga_p<2<r\ga_r\le q\ga_q$ and again proceeding as {\it Case 1}, by combining $W(r\ga_r)=0$ and $W(p\ga_p)=0$, we can
	obtain the constant $\al_2$ with the required properties. If $r>q$, then there is $p\ga_p<2<q\ga_q<r\ga_r$ and
	again proceeding as {\it Case 1}, by combining $W(q\ga_q)=0$ and $W(p\ga_p)=0$, we can obtain the constant $\al_2$ with the
	required properties.
\end{itemize}
\ep

	Using Lemmas \ref{201905-4h}, \ref{201905-4pho}, we can discribe the geometry of $\p_{a_1,a_2}$.
\bl\lab{201905-4str}
	If $T(a_1,a_2)<\min\{\al_1,\al_2\}$, then for every $(u,v)\in S_{a_1}\times S_{a_2}$, the function $\Phi_{(u,v)}(t)$
	has exactly two critical points $s_{(u,v)}<t_{(u,v)}$ and two zero-points $c_{(u,v)}<d_{(u,v)}$
	with $s_{(u,v)}<c_{(u,v)}<t_{(u,v)}<d_{(u,v)}$. Moreover:
	\begin{itemize}[fullwidth,itemindent=0em]
		\item[(1)]$s\s(u,v)\in \p^+_{a_1,a_2}$ iff $s=s_{(u,v)}$; $s\s(u,v)\in\p^-_{a_1,a_2}$ iff $s=t_{(u,v)}$;
		\item[(2)]$s_{(u,v)}<\log\f{R_0}{(\int_{\RN}|\nabla u|^2+|\nabla v|^2)^{1/2}}$ and
			$$\Phi_{(u,v)}(s_{(u,v)})=\inf\ \left\{\Phi_{(u,v)}(s):s\in\big(-\iy,\log\f{R_0}{(\int_{\RN}|\nabla u|^2+|\nabla v|^2)^{1/2}}\big)\right\};$$
		\item[(3)]$I\big(t_{(u,v)}\s(u,v)\big)=\max_{s\in\R}I\big(s\s(u,v)\big)>0$;
		\item[(4)]the maps $(u,v)\mapsto t_{(u,v)}$ and $(u,v)\mapsto s_{(u,v)}$ are of class $C^1$.
	\end{itemize}
\el
\bp
Let $(u,v)\in S_{a_1}\times S_{a_2}$. By \eqref{201911-11}, we have
	$$\Phi_{(u,v)}(s)=I\big(s\s(u,v)\big)\geq h\big(e^s(\int_{\RN}|\nabla u|^2+|\nabla v|^2)^{1/2}\big),$$
and then 
	$$\Phi_{(u,v)}(s)>0,\quad \forall s\in\big(\log\f{R_0}{(\int_{\RN}|\nabla u|^2+|\nabla v|^2)^{1/2}},
		\log\f{R_1}{(\int_{\RN}|\nabla u|^2+|\nabla v|^2)^{1/2}}\big).$$
Recalling the facts that $\Phi_{(u,v)}(-\infty)=0^-$ and $\Phi_{(u,v)}(+\infty)=-\infty$,
we see that $\Phi_{(u,v)}$ has at least two critical points $s_{(u,v)}<t_{(u,v)}$,  where $s_{(u,v)}$ is local minimum point
on $\big(-\iy,\log\f{R_0}{(\int_{\RN}|\nabla u|^2+|\nabla v|^2)^{1/2}}\big)$ at negetive level, and $t_{(u,v)}$ ia a global maximum point
at positive level. On the other hand, from Remark \ref{Remark1}, $\Phi_{(u,v)}(s)$ has at most two critical points in $(-\iy,+\iy)$, which means
that $\Phi_{(u,v)}(s)$ has exactly two critical points $s_{(u,v)}$ and $t_{(u,v)}$. Since $\Phi_{(u,v)}'(s)=P\big(s\s(u,v)\big)$, we
know that  $s\s(u,v)\in\p_{a_1,a_2}$ implies $s=s_{(u,v)}\ \text{or}\ t_{(u,v)}$. Moreover, noticing $\Phi_{(u,v)}''(s_{(u,v)})\geq0$,
 $\Phi_{(u,v)}''(t_{(u,v)})\leq0$ and $\p^0_{a_1,a_2}=\emptyset$, we deduce that $s_{(u,v)}\s(u,v)\in\p^+_{a_1,a_2}$ and
  $t_{(u,v)}\s(u,v)\in\p^-_{a_1,a_2}$.

\vskip0.1in
	By the monotonicity, $\Phi_{(u,v)}$ has exactly two zero-points $c_{(u,v)}$ and $d_{(u,v)}$, with $s_{(u,v)}<c_{(u,v)}<t_{(u,v)}<d_{(u,v)}$.
It remains to show that the maps $(u,v)\mapsto t_{(u,v)}$ and $(u,v)\mapsto s_{(u,v)}$ are of class $C^1$. We apply the implicit function
theorem on $\Psi(s,u,v)=\Phi_{(u,v)}'(s)$. Using the fact that
	$$\Psi(s_{(u,v)},u,v)=\Psi(t_{(u,v)}\s(u,v))=0,$$
	$$\partial_s\Psi(s_{(u,v)},u,v)=\Phi_{(u,v)}''(s_{(u,v)})>0,$$
	$$\partial_s\Psi(s_{(u,v)},u,v)=\Phi_{(u,v)}''(t_{(u,v)})<0,$$
and $\p^0_{a_1,a_2}=\emptyset$, there is $(u,v)\mapsto t_{(u,v)}$ and $(u,v)\mapsto s_{(u,v)}$ are of class $C^1$.
\ep

	For $k>0$, let
	$$A_R:=\big\{(u,v)\in S_{a_1}\times S_{a_2}:(\int_{\RN}|\nabla u|^2+|\nabla v|^2)^{1/2}<R\big\}.$$
We have the following crucial estimates.

\bl\lab{201906-4upbdd}
	Let $N\ge3$ and $T(a_1,a_2)<\min\{\al_1,\al_2\}$. If $r_2<2$, then
		$$m(a_1,a_2)=\inf_{(u,v)\in A_{R_0}}I(u,v)<\min\big\{m(a_1,0),m(0,a_2)\big\}.$$
\el
\bp
From Lemma \ref{201905-4str}, we have
	$$\p^+_{a_1,a_2}=\left\{s_{(u,v)}\star(u,v):(u,v)\in S_{a_1}\times S_{a_2}\right\}\subset A_{R_0},$$
and
	$$m(a_1,a_2)=\inf_{\p_{a_1,a_2}} I(u,v)=\inf_{\p^+_{a_1,a_2}} I(u,v)<0.$$
Obviously $m(a_1,a_2)\ge\inf_{A_{R_0}}I(u,v)$. On the other hand, for any $(u,v)\in A_{R_0}$,  since
 $0<\log \f{R_0}{(\int_{\RN}|\nabla u|^2+|\nabla v|^2)^{1/2}}$, there is 
	$$m(a_1,a_2)\leq I(s_{(u,v)}\star(u,v))\le I(u,v).$$
It follows that $m(a_1,a_2)=\inf_{(u,v)\in A_{R_0}}I(u,v)$. Noting that $p<\bar p<q$ means $m(a_1,0)<0<m(0,a_2)$, we only need to 
prove $m(a_1,a_2)<m(a_1,0)$.

	We now choose a proper test function to prove $m(a_1,a_2)<m(a_1,0)$. From $h(R_0)=0$, we have $\f{1}{2}R_0^2>\D_2R_0^{p\ga_p}$,
that is $R_0^{2-p\ga_p}>2\D_2$. Let $(u,\la_0)\in S_{a_1}\times \R^+$ be the unique function in Lemma \ref{201905-nlses} with parameters $p,\mu_1,a_1$. It follows that
	$$|\nabla u|_2^2=\ga_p\mu_1|u|_p^p\le p\ga_p\D_2|\nabla u|_2^{p\ga_p}<R_0^{2-p\ga_p}|\nabla u|_2^{p\ga_p},$$
which means $|\nabla u|_2<R_0$. Take $m>1$ such that $\f{N}{2}-\f{2}{r_2}<m<\f{N}{2}-1$ and let 
	$$\vp(x)\in C_0^\iy(B_2(0)),\quad 0\le\vp(x)\le1,\quad\vp(x)=1\ \text{in}\ B_1(0).$$
We define $v(x)=c\f{\vp(x)}{|x|^m}\in H^1$ with constant $c>0$ such that $v\in S_{a_2}$.
Therefore, $(u,s\s v)\in A_{R_0}$ for $s\ll-1$. Let
	$$\al(s)=\int_{\RN}|u|^{r_1}|s\s v|^{r_2}=C_0e^{(\f{N}{2}-m)r_2s}\int_{\RN}u^{r_1}(x)\f{\vp^{r_2}(e^sx)}{|x|^{mr_2}}.$$
From the decay properties of \cite{Ni=CPDE=1993}, we get by replacing a variable that $u$ decays exponentially
	$$u(x)=O(|x|^{-\frac{N-2}{2}}e^{-\la_0^{\f{1}{2}}|x|}),\quad \text{as}\ |x|\ra\iy,$$
and $|u(x)|\le M$ in $\RN$. Then
	$$0<\int_{\RN}\f{u^{r_1}(x)}{|x|^{mr_2}}\le C\Big(\int_{B_R(0)}\f{1}{|x|^{mr_2}}+
		\int_{B_R(0)^c}|x|^{-\frac{(N-2)r_1}{2}-mr_2}e^{-\la_0^{\f{1}{2}}r_1|x|}\Big)<\iy.$$
Thus by the Dominated Convergence Theorem, we obtain
	$$\lim_{s\to-\iy}\int_{\RN}u^{r_1}(x)\f{\vp^{r_2}(e^sx)}{|x|^{mr_2}}=\int_{\RN}u^{r_1}(x)\f{1}{|x|^{mr_2}}=:C_1\in(0,\iy).$$
Hence $\al(s)=C_0e^{\theta s}(C_1+o(1))>\f{C_0C_1}{2}e^{\theta s}$ as $s\ra-\iy$ where $\theta=(\f{N}{2}-m)r_2\in(1,2)$. 
Finally we see that for some $s\ll-1$, there  holds
\bg
	m(a_1,a_2)&\le I(u,s\s v)\\
		&=E_{p,a_1,\mu}(u)+\f{e^{2s}}{2}|\nabla v|_2^2-\f{e^{q\ga_qs}}{q}|v|_q^q-\beta\al(s)\\
		&<E_{p,a_1,\mu}(u)+\f{e^{2s}}{2}|\nabla v|_2^2-\f{e^{q\ga_qs}}{q}|v|_q^q-\beta\f{C_0C_1}{2}e^{\theta s}\\
		&<E_{p,a_1,\mu}(u)=m(a_1,0).
\end{align*}
\ep

	Now we prove the compactness of Palais-Smale sequences.

\bl\lab{201906-4cpt}
	Suppose $N=3,4$ and $\D_1\beta<\f{1}{4}$ when $r=\bar p$. Let $(u_n,v_n)\subset S_{a_1}\times S_{a_2}$ be  a radial
	Palais-Smale sequence for $I|_{S_{a_1}\times S_{a_2}}$ at level $m(a_1,a_2)$ with additional properties $P(u_n,v_n)\ra0$
	and $u_n^-,v_n^-\ra0$ a.e. in $\RN$. If
	$$m(a_1,a_2)<\min\big\{m(a_1,0),m(0,a_2)\big\},$$
	then up to a subsequence, $(u_n,v_n)\to(u,v)$ in $H$, where $(u,v)$ is a positive solution of \eqref{201905-2} for
	some $\la_1,\la_2>0$.
\el
\bp
	We first prove that $(u_n,v_n)$ is bounded. Let $\rho_n=(|u_n|_2^2+|v_n|_2^2)^{\f{1}{2}}$ and
	$$\begin{aligned}
		Z_n(t):&=tI(u_n,v_n)-P(u_n,v_n)\\
		&=\f{t-2}{2}\int|\nabla u_n|^2+|\nabla v_n|^2-\f{t-p\ga_p}{p}\int\mu_1 |u_n|^p\\
		&\quad-\f{t-q\ga_q}{q}\int\mu_2|v_n|^q-(t-r\ga_r)\int\beta|u_n|^{r_1}|v_n|^{r_2}\\
		&\le C(t), \quad\quad \forall\ n\ge1.
	  \end{aligned} $$

	We still disscuss it in four different situations.
\begin{itemize}[fullwidth,itemindent=0em]
\item[{\it Case 1: $r<\bar p$.}]
	From $Z_n(q\ga_q)\le C$, we get
		$$\f{q\ga_q-2}{2}\rho_n^2\le C+\f{q\ga_q-p\ga_p}{p}\int\mu_1 |u_n|^p+(q\ga_q-r\ga_r)\int\beta|u_n|^{r_1}|v_n|^{r_2}\le C(1+\rho_n^{p\ga_p}+\rho_n^{r\ga_r}),$$
	which implies that $(u_n,v_n)$ is bounded.
\item[{\it Case 2: $r=\bar p$.}]
	Note that $r\ga_r=2$. From $Z_n(q\ga_q)\le C$, we get
	\bg
		\f{q\ga_q-2}{2}(1-2\D_1\beta)\rho_n^2&\le C+\f{q\ga_q-2}{2}\rho_n^2-(q\ga_q-2)\int\beta|u_n|^{r_1}|v_n|^{r_2}\\
		&\le C+\f{q\ga_q-p\ga_p}{p}\int\mu_1 |u_n|^p\le C(1+\rho_n^{p\ga_p}),
	\end{align*}
	which implies that $(u_n,v_n)$ is bounded.
\item[{\it Case-3: $\bar p<r\le q$.}]
	From $Z_n(r\ga_r)\le C$, we get
		$$\f{r\ga_r-2}{2}\rho_n^2\le C+\f{r\ga_r-p\ga_p}{p}\int\mu_1 |u_n|^p\le C(1+\rho_n^{p\ga_p}),$$
	which implies that $(u_n,v_n)$ is bounded.
\item[{\it Case-4: $\bar p<q<r$.}]	From $Z_n(q\ga_q)\le C$, we get
		$$\f{q\ga_q-2}{2}\rho_n^2\le C+\f{q\ga_q-p\ga_p}{p}\int\mu_1 |u_n|^p\le C(1+\rho_n^{p\ga_p}),$$
	which implies that $(u_n,v_n)$ is bounded.
\end{itemize}

\vskip0.08in
	Since the sequence $(u_n,v_n)$ is a bounded sequence of radial functions, by the compactness of the embedding $H^1_r\hra L^p(\RN)$
for $2<p<2^*$, there exists a $(u,v)\in H$ such that up to a subsequence $(u_n,v_n)\rh(u,v)$ in $H$ and $L^{2^*}(\RN)\times L^{2^*}(\RN)$  and $(u_n,v_n)\ra(u,v)$
in $L^p(\RN)\times L^p(\RN),L^r(\RN)\times L^r(\RN),L^q(\RN)\times L^q(\RN)$ when $q<2^*$, and $(u_n,v_n)\to(u,v)$ a.e. in $\RN$. Hence $u,v\geq0$ are
radial funtions. Since $I|_{S_{a_1}\times{S_{a_2}}}'(u_n,v_n)\ra 0$, by the Lagrange multiplier's rule, we have that there exists a
sequence $(\la_{1,n},\la_{2,n})\subset \R^2$ such that
\be\lab{201906-8.1}
	\int_{\RN}\nabla u_n\cdot\nabla\vp+\la_{1,n}u_n\vp-\mu_1|u_n|^{p-2}u_n\vp-\beta r_1|u_n|^{r_1-2}|v_n|^{r_2}u_n\vp=o(1)||\varphi||_{H^1},
\ee
\be\lab{201906-8.2}
	\int_{\RN}\nabla v_n\cdot\nabla \psi+\la_{2,n}v_n\psi-\mu_2|v_n|^{q-2}v_n\psi-\beta r_2|u_n|^{r_1}|v_n|^{r_2-2}v_n\psi=o(1)||\psi||_{H^1},
\ee
as $n\to\iy$, for every $(\vp,\psi)\in H$. By choosing $(\vp,\psi)=(u_n,v_n)$, we decude that $(\lambda_{1,n},\lambda_{2,n})$ is
bounded as well, and hence up to a subsequence $(\lambda_{1,n},\lambda_{2,n})\ra(\lambda_1,\lambda_2)\in \R^2$. Then, passing
to the limits in \eqref{201906-8.1}-\eqref{201906-8.2}, we deduce that $(u,v)$ is a nonnegative solution of \eqref{201905-2}.
Thus from the Pohozaev identity we obtain
\be\lab{201906-9}
\la_1|u|_2^2+\la_2|v|_2^2=(1-\ga_p)\int_{\RN}\mu_1u^p+(1-\ga_q)\int_{\RN}\mu_2v^q+(1-\ga_r)r\int_{\RN}\beta u^{r_1}v^{r_2}.
\ee
Moreover, combining $P(u_n,v_n)\ra0$ with \eqref{201906-8.1}\eqref{201906-8.2}, we have
\be\lab{201906-10}
	\begin{aligned}
	&\quad\ \la_1a_1^2+\la_2a_2^2=\lim_{n\to\iy}\la_{1,n}|u_n|_2^2+\la_{2,n}|v_n|_2^2\\
	&=\lim_{n\to\iy}\int_{\RN}-(|\nabla u_n|^2+|\nabla v_n|^2)+\mu_1 |u_n|^p+\mu_2 |v_n|^q+r\beta |u_n|^{r_1}|v_n|^{r_2}\\
	&=\lim_{n\to\iy}(1-\ga_p)\int_{\RN}\mu_1 |u_n|^p+(1-\ga_q)\int_{\RN}\mu_2 |v_n|^q+(1-\ga_r)r\int_{\RN}\beta |u_n|^{r_1}|v_n|^{r_2}\\
	&=(1-\ga_p)\int_{\RN}\mu_1u^p+(1-\ga_q)\int_{\RN}\mu_2v^q+(1-\ga_r)r\int_{\RN}\beta u^{r_1}v^{r_2}.
	\end{aligned}
\ee

	Now since there may be $u=0$ or $v=0$, we will proceed in four cases.
\begin{itemize}[fullwidth,itemindent=0em]
\item[{\it Case 1: $u=0,v=0$.}]
	Since $(u_n,v_n)\ra(u,v)$ in $L^p(\RN)\times L^p(\RN),L^r(\RN)\times L^r(\RN)$, we have
		$$0=P(u_n,v_n)+o(1)=\int_{\RN}|\nabla u_n|^2+|\nabla v_n|^2-\ga_q\int_{\RN}\mu_2|v_n|^q+o(1).$$
	Then it follows that
	\bg
		m(a_1,a_2)&=\lim_{n\ra\iy}I(u_n,v_n)=\lim_{n\ra\iy}\f{1}{2}\int_{\RN}|\nabla u_n|^2+|\nabla v_n|^2-\f{1}{q}\int_{\RN}\mu_2|v_n|^q\\
		&=\lim_{n\ra\iy}(\f{1}{2}-\f{1}{q\ga_q})\int_{\RN}|\nabla u_n|^2+|\nabla v_n|^2\ge0.
	\end{align*}
	However, $m(a_1,a_2)<m(a_1,0)<0$, we get a contradiction.

\vskip0.08in
\item[{\it Case 2: $u\ne0,v=0$.}]
	By the maximum principle, $u$ is a positive solution of \eqref{201905-nlse} with parameters $p,\mu_1$ and $a=|u|_2\le a_1$, 
	then $m(a_1,0)\le m(|u|_2,0)=I(u,0)$. Let $\bar u_n=u_n-u$, then by  using the Brezis-Lieb Lemma and Lemma \ref{201905-BL}, we have
	\bg
		0&=P(u_n,v_n)+o(1)=P(\bar u_n,v_n)+P(u,0)+o(1)\\
		&=\int_{\RN}|\nabla \bar u_n|^2+|\nabla v_n|^2-\ga_q\int_{\RN}\mu_2|v_n|^q+o(1),
	\end{align*}
	and hence
	\bg
		m(a_1,a_2)&=\lim_{n\ra\iy}I(u_n,v_n)=\lim_{n\ra\iy}I(\bar u_n,v_n)+I(u,0)\\
		&\ge\lim_{n\ra\iy}\f{1}{2}\int_{\RN}|\nabla\bar u_n|^2+|\nabla v_n|^2-\f{1}{q}\int_{\RN}\mu_2|v_n|^q+m(a_1,0)\\
		&=\lim_{n\ra\iy}(\f{1}{2}-\f{1}{q\ga_q})\int_{\RN}|\nabla\bar u_n|^2+|\nabla v_n|^2+m(a_1,0)\ge m(a_1,0),
	\end{align*}
	which is a contradiction.

\vskip0.08in
\item[{\textit Case 3: $u=0,v\ne0$.}]
	If $q=2^*$, then $v$ is a positive solution of \eqref{201905-nlse} with parameters $p=2^*,\mu=\mu_2$ and $a=|v|_2>0$, which 
	contradicts Lemma \ref{201905-nlsec}. If $q<2^*$, then proceeding as {\it Case 2}, we get $m(a_1,a_2)\ge m(0,a_2)$, a contradiction too.

\vskip0.08in
\item[{\textit Case 4: $u\ne0,v\ne0$.}]
	In this case, we prove $(u_n,v_n)\ra(u,v)$ in $H$. Again by the maximum principle, $u,v>0$, then Lemma \ref{201905-liouville}
	implies $\la_1,\la_2>0$. Moreover, from \eqref{201906-9}\eqref{201906-10}, we obtain
		$$\la_1(a_1^2-|u|_2^2)+\la_2(a_2^2-|v|_2^2)=0,$$
	and since $0<|u|_2\le a_1,0<|v|_2\le a_2$ there must be $|u|_2=a_1,|v|_2=a_2$. So $(u,v)\in\p_{a_1,a_2}$. Let $(\bar u_n,\bar v_n)=(u_n-u,v_n-v)$, then we have
	\bg
		0&=P(u_n,v_n)+o(1)=P(\bar u_n,\bar v_n)+P(u,v)+o(1)\\
		&=\int_{\RN}|\nabla\bar u_n|^2+|\nabla\bar v_n|^2-\ga_q\int_{\RN}\mu_2|\bar v_n|^q+o(1),
	\end{align*}
	and hence
	\bg
	m(a_1,a_2)&=\lim_{n\ra\iy}I(u_n,v_n)=\lim_{n\ra\iy}I(\bar u_n,\bar v_n)+I(u,v)\\
		&\ge\lim_{n\ra\iy}\f{1}{2}\int_{\RN}|\nabla\bar u_n|^2+|\nabla\bar v_n|^2-\f{1}{q}\int_{\RN}\mu_2|\bar v_n|^q+m(a_1,a_2)\\
		&=\lim_{n\ra\iy}(\f{1}{2}-\f{1}{q\ga_q})\int_{\RN}|\nabla\bar u_n|^2+|\nabla\bar v_n|^2+m(a_1,a_2)\ge m(a_1,a_2).
	\end{align*}
	So $I(u,v)=m(a_1,a_2)$ and $(u_n,v_n)\ra(u,v)$ in $H$.
\end{itemize}
\ep

\vskip0.12in
\bp[Proof of the Theorem \ref{thm2}]
	Take $\al_0=\min\{\al_1,\al_2\}$, then by Lemma \ref{201906-4upbdd} and Lemma \ref{201906-4cpt}, it is sufficient to prove the existence
of a radial Palais-Smale sequence for $I|_{S_{a_1}\times S_{a_2}}$ at level $m(a_1,a_2)$ with additional properties $P(u_n,v_n)\ra0$
and $u_n^-,v_n^-\ra0$ a.e. in $\RN$.

	Let $m_r(a_1,a_2)=\inf_{A_{R_0}\cap H_r} I(u,v)$, and using the symmetric decreasing rearrangement it is easy to check
 $m(a_1,a_2)=m_r(a_1,a_2)$. Choosing a minimizing sequence $(\tilde u_n,\tilde v_n)$ for $m(a_1,a_2)=\inf_{A_{R_0}\cap H_r} I(u,v)$,
we assume $(\tilde u_n,\tilde v_n)$ are nonnegative by replacing  $(\tilde u_n,\tilde v_n)$ with $(|\ti u_n|,|\ti v_n|)$.
Futhermore, using the fact that $I\big(s_{(\tilde u_n,\tilde v_n)}\s(\tilde u_n,\tilde v_n)\big)\leq I(\tilde u_n,\tilde v_n)$, 
and replacing $(\tilde u_n,\tilde v_n)$ by $s_{(\tilde u_n,\tilde v_n)}\star(\tilde u_n,\tilde v_n)$, we obtain a minimizing sequence
$(\tilde u_n,\tilde v_n)\in\p^+_{a_1,a_2,r}$ for $n\ge1$. Therefore, by Ekeland's varational principle, there is a radial
Palais-Smale sequence $(u_n,v_n)$ for $I|_{S_{a_1,r}\times S_{a_2,r}}$ (hence a Palais-Smale sequence for
 $I|_{S_{a_1}\times S_{a_2}}$) with the property $||(u_n,v_n)-(\tilde u_n,\tilde v_n)||\to0$ as $n\to\infty$, which implies that
	$$P(u_n,v_n)=P(\tilde u_n,\tilde v_n)+o(1)\to0\quad \text{and}\quad u_n^-,v_n^-\ra0\ \text{a.e. in}\ \RN,$$
then we finish the proof.
\ep

\vskip0.23in
\section{The purely $L_2$-supercritical case}

	In this section, we suppose $\bar p<p,q,r<2^*$. To start our discussion, we consider once again the Pohozaev manifold $\p_{a_1,a_2}$
and the decomposition $\p_{a_1,a_2}=\p^+_{a_1,a_2}\cup\p_{a_1,a_2}^0\cup\p_{a_1,a_2}^-$. If there is a $(u,v)\in\p_{a_1,a_2}^0$,
then combining $\Phi_{(u,v)}'(0)=0$ and $\Phi_{(u,v)}''(0)=0$, we deduce that
	$$(p\ga_p-2)\ga_p\int_{\RN}\mu_1 |u|^p+(q\ga_q-2)\int_{\RN}\mu_2|v|^q+(r\ga_r-2)r\ga_r\int_{\RN}\beta |u|^{r_1}|v|^{r_2}=0.$$
Since $p\ga_p,r\ga_r,q\ga_q>2$, there must be $(u,v)=(0,0)$, in contradiction with $(u,v)\in S_{a_1}\times S_{a_2}$. This shows
that $\p_{a_1,a_2}^0=\emptyset$, and then as Lemma \ref{201905-4pho} we can prove that $\p_{a_1,a_2}$ is a $C^1$ submanifold in $H$ with codimension 3.
However, in this section, the geometry of $\p_{a_1,a_2}$ will be different from the one in Lemma \ref{201905-4str}.

\bl\lab{201906-6str}
	For any $(u,v)\in S_{a_1}\times S_{a_2}$, the function $\Phi_{(u,v)}$ has a unique critical point $t_{(u,v)}\in\R$, which is a
strict maximum point at positive level. Moreover,
\begin{itemize}[fullwidth,itemindent=0em]
	\item[(1)] $\p_{a_1,a_2}=\p_{a_1,a_2}^-$ and $P(u,v)<0$ iff $t_{(u,v)}<0$;
	\item[(2)] $\Phi_{(u,v)}$ is strict increasing in $(-\iy,t_{(u,v)})$;
	\item[(3)] the map $(u,v)\mapsto t_{(u,v)}$ is of class $C^1$.
\end{itemize}
\el
\bp
The proof is completely the sames as \cite[Lemma 6.1]{Soave=arXiv=2019}, and we omit the details here.
\ep

	Using the above lemma, it is easy to see that
	$$m(a_1,a_2)=\inf_{S_{a_1}\times S_{a_2}}\max_{t\in\R}I(t\s(u,v)).$$
By using  the same techniques as that in  Lemma \ref{201906-4cpt}, we can prove the following lemma.

\bl\lab{201906-6cpt}
	Suppose $N=3,4$. Let $(u_n,v_n)\subset S_{a_1}\times S_{a_2}$ be a radial Palais-Smale sequence for $I|_{S_{a_1}\times S_{a_2}}$
	at level $m(a_1,a_2)$ with the additional properties $P(u_n,v_n)\ra0$ and $u_n^-,v_n^-\ra0$ a.e. in $\RN$. If
	$$0<m(a_1,a_2)<\min\big\{m(a_1,0),m(0,a_2)\big\},$$
	then up to a subsequence $(u_n,v_n)\to(u,v)$ in $H$, where $(u,v)$ is a positive solution of \eqref{201905-2} for some $\la_1,\la_2>0$.
\el
\br
It is natural  that $m(a_1,a_2)>0$. Indeed, for any $(u,v)\in\p_{a_1,a_2}$, there is
\bg
	\int_{\RN}|\nabla u|^2+|\nabla v|^2&=\ga_p\int_{\RN}\mu_1|u|^p+\ga_q\int_{\RN}\mu_2|v|^q+r\ga_r\int_{\RN}\beta |u|^{r_1}|v|^{r_2}\\
		&\le\D_2p\ga_p(\int_{\RN}|\nabla u|^2)^{\f{p\ga_p}{2}}+\D_3q\ga_q(\int_{\RN}|\nabla v|^2)^{\f{q\ga_q}{2}}\\
		&\quad\ +\D_1r\ga_r\beta(\int_{\RN}|\nabla u|^2+|\nabla v|^2)^{\f{r\ga_r}{2}},
\end{align*}
which implies $\inf_{\p_{a_1,a_2}}\int|\nabla u|^2+|\nabla v|^2\ge C>0$. So we have
\bg
	&\quad m(a_1,a_2)=\inf_{\p_{a_1,a_2}}I(u,v)\\
	&=\inf_{\p_{a_1,a_2}}\f{p\ga_p-2}{2p}\int_{\RN}\mu_1 |u|^p+\f{q\ga_q-2}{2q}\int_{\RN}\mu_2 |u|^q
			+\f{r\ga_r-2}{2}\int_{\RN}\beta |u|^{r_1}|v|^{r_2}\\
	&\ge C\inf_{\p_{a_1,a_2}}\int_{\RN}|\nabla u|^2+|\nabla v|^2>0.
\end{align*}
\er

	We recall the following lemma in \cite{BartschSoave=2017}.
\bl\lab{201906-con}
	The map $(s,u)\in\R\times H^1\to s\s u\in H^1$ is continuous.
\el
	Now we give a way to find such a Palais-Smale sequence as the required one in Lemma \ref{201906-6cpt}.
\bl\lab{201909-6PSseq}
	There is a radial Palais-Smale sequence for $I|_{S_{a_1}\times S_{a_2}}$ at level $m(a_1,a_2)$ with the additional properties
	 $P(u_n,v_n)\ra0$ and $u_n^-,v_n^-\ra0$ a.e. in $\RN$.
\el
\bp
We consider the functional $\ti I:\R\times H^1(\RN)\times H^1(\RN)\to\R$ defined by
	$$\ti I(s,u,v):=I(s\s(u,v))$$
on the constraint $\R\times S_{a_1,r}\times S_{a_2,r}$. Denote the closed sublevel set
by $I^c=\{(u,v)\in S_{a_1}\times S_{a_2}:I(u,v)\leq c\}$. We note that for any $(u,v)\in S_{a_1}\times S_{a_2}$,
	$$I(u,v)\geq\f{1}{2}(|\nabla u|_2^2+|\nabla v|_2^2)-\D_2|\nabla u|_2^{p\ga_p}-\D_3|\nabla v|_2^{q\ga_q}-
	\D_1\beta(|\nabla u|_2^2+|\nabla v|_2^2)^{\f{r\ga_r}{2}},$$
	$$I(u,v)\leq\frac{1}{2}(|\nabla u|_2^2+|\nabla v|_2^2),$$
	$$P(u,v)\geq|\nabla u|_2^2+|\nabla v|_2^2-\D_1p\ga_p|\nabla u|_2^{p\ga_p}-\D_3q\ga_q|\nabla v|_2^{q\ga_q}-
	\D_1r\ga_r\beta(|\nabla u|_2^2+|\nabla v|_2^2)^{\f{r\ga_r}{2}},$$
then there exists a small $k>0$ such that
	$$0<I(u,v)<m(a_1,a_2),\quad P(u,v)>0,\quad \forall (u,v)\in \bar A_{k}.$$
We introduce the minimax class
	$$\Gamma:=\{\ga=(\al,\vp_1,\vp_2)\in C([0,1],\R\times S_{a_1,r}\times S_{a_2,r}):\ga(0)\in\{0\}\times\bar A_k,\ga(1)\in\{0\}\times I^0\}$$
with the associated minimax level
	$$\sigma:=\inf_{\ga\in\Gamma}\max_{t\in[0,1]}\ti I(\ga(t)).$$
Next we check that $\sigma=m(a_1,a_2)$. On the one hand, for any $(u,v)\in\p_{a_1,a_2}$, there are $(u^*,v^*)\in S_{a_1,r}\times S_{a_2,r}$ and
 $P(u^*,v^*)\leq P(u,v)=0$, which implies $t_*=t_{(u^*,v^*)}\leq0$. It follows that
	$$I(u,v)\ge I(t_*\s(u,v))\ge I(t_*\s(u^*,v^*))=\max_{t\in\R} I(t\s(u^*,v^*)).$$
Observing that
	$$|\nabla s\s u^*|_2^2+|\nabla s\s v^*|_2^2\to0,\quad \text{as}\ s\ra-\iy,$$
	$$I(t\s(u^*,v^*))\ra-\iy, \quad \text{as}\ s\ra\iy,$$
we choose $s_0\ll-1,s_1\gg1$ such that $s_0\s(u^*,v^*)\in A_k$ and $s_1\s(u^*,v^*)\in I^0$.
Then we define $\ga_*:[0,1]\to\R\times S_{a_1,r}\times S_{a_2,r}$ 
	$$\ga_*(t)=\big(0,[(1-t)s_0+ts_1]\star(u^*,v^*)\big),$$
and by Lemma \ref{201906-con}, $\ga_*\in\Gamma$. Hence
	$$\sigma\leq\max_{t\in[0,1]}\ti I(\ga_*(t))\leq\max_{t\in\R}I(t\s(u^*,v^*))\leq I(u,v),$$
which implies $\sigma\leq m(a_1,a_2)$. On the other hand, for any $\ga=(\al,\vp_1,\vp_2)\in\Gamma$, we consider the function
	$$P_\ga:t\in[0,1]\to P\big(\al(t)\s(\vp_1(t),\vp_2(t))\big)\in\R.$$
It is easy to see that $P_\ga$ is continuous and $P_\ga(0)>0$. We claim that $P_\ga(1)<0$. Indeed, if $P_\ga(1)\geq0$,
we have $t_{(\vp_1(1),\vp_2(1))}\ge0$, and then from Lemma \ref{201906-6str},
	$$I(\vp_1(1),\vp_2(1))=\Phi_{(\vp_1(1),\vp_2(1))}(0)>\Phi_{(\vp_1(1),\vp_2(1))}(-\iy)=0^+,$$
which is a contradiction. Thus we obtain a $t_\ga\in(0,1)$ such that $P_\ga(t_\ga)=0$. It follows that
	$$\max_{t\in[0,1]}\ti I(\ga(t))\ge\ti I(\ga(t_\ga))=I(\al(t_\ga)\s(\vp_1(t_\ga),\vp_2(t_\ga)))\ge m(a_1,a_2),$$
which implies $\sigma\geq m(a_1,a_2)$. Hence $\sigma=m(a_1,a_2)$.

	Let $\mc F=\{\ga([0,1]):\ga\in\Gamma\}$. Using the terminology in \cite[Section 5]{Ghoussoub=1993}, $\mc F$ is a homotopy stable
family of compact subset of $\R\times S_{a_1,r}\times S_{a_2,r}$ with extended closed boundary $\{0\}\times \bar A_k\cup \{0\}\times I^0$,
and the superlevel set $\{\ti I\ge\sigma\}$ is a dual set for $\mc F$, which means that the assumptions in \cite[Theorem 5.2]{Ghoussoub=1993}
are satisfied. Therefore, taking a minimizing sequence $\left\{\ga_n([0,1]),\ga_n=(\al_n,\vp_{1,n},\vp_{2,n})\right\}$ for $\sigma$
with the property that $\al(t)=0$, $\vp_{1,n}(t)\geq0$, $\vp_{2,n}(t)\geq0$ for every $t\in[0,1]$( Indeed, we can replace $\ga_n$
by $\ti \ga_n=(0,\al_n\s(|\vp_{1,n}|,|\vp_{2,n}|))$ ), there exists a sequence $(s_n,u_n,v_n)\subset\R\times S_{a_1,r}\times S_{a_2,r}$
such that as $n\ra\iy$, $\ti I(s_n,u_n,v_n)\ra\sigma$ and
\be\lab{201906-14.1}
	\pa_s\ti I(s_n,u_n,v_n)\ra0,\quad ||\pa_{(u,v)}\ti I(s_n,u_n,v_n)||_{T_{u_n}S_{a_1,r}\times T_{v_n}S_{a_2,r}}\ra0,
\ee
\be\lab{201906-14.2}
	|s_n|+\text{dist}\left((u_n,v_n),(\vp_{1,n}([0,1]),\vp_{2,n}([0,1]))\right)\ra0.
\ee
Let $(\bar u_n,\bar v_n)=s_n\s(u_n,v_n)\in S_{a_1,r}\times S_{a_2,r}$. From \eqref{201906-14.2}, we know that $\{s_n\}$ is bounded
and $\bar u_n^-,\bar v_n^-\to0$ a.e. in $\RN$. Moreover, \eqref{201906-14.1} implies that
	$$P(\bar u_n,\bar v_n)=\partial_s\tilde I(s_n,u_n,v_n)\to0,$$
and that
\bg
	I'(\bar u_n,\bar v_n)[\phi,\psi]&=\pa_{(u,v)}\ti I(s_n,u_n,v_n)[(-s_n)\s(\phi,\psi)]\\
		&=o(1)||(-s_n)\s(\phi,\psi)||_H\\
		&=o(1)||(\phi,\psi)||_H,
\end{align*}
for any $(\phi,\psi)\in T_{\bar u_n}S_{a_1,r}\times T_{\bar v_n}S_{a_2,r}$.
Summing up, $(\bar u_n,\bar v_n)$ is a radial Palais-Smale sequence of $I|_{S_{a_1}^r\times S_{a_2}^r}$ and hence a radial symmetric
Palais-Smale sequence of $I|_{S_{a_1}\times S_{a_2}}$ at level $\sigma$.
\ep

	Before giving an estimate of $m(a_1,a_2)$ from above, we would like to study the dependence
of $m(a_1,a_2)$ on $\beta$. In the following lemma, we denote $m(a_1,a_2),I(u,v)$ by $m_\beta(a_1,a_2)$, $I_\beta(u,v)$ respectively.
\bl\lab{201909-6mbeta}
	For any $a_1,a_2>0$,  we have that 
	\begin{itemize}[fullwidth,itemindent=0em]
	\item[(1)]$m_\beta(a_1,a_2)$ is decreasing with respect to $\beta\ge0$;
	\item[(2)]$m_0(a_1,a_2)=\min\big\{m(a_1,0),m(0,a_2)\big\}$.
	\end{itemize}
\el
\bp
\begin{itemize}[fullwidth,itemindent=0em]
\item[(1)] For any  $ \beta_1\ge\beta_2\ge0$,
		$$m_{\beta_1}(a_1,a_2)=\inf_{S_{a_1}\times S_{a_2}}\max_{t\in\R}I_{\beta_1}(t\s(u,v))\le \inf_{S_{a_1}\times S_{a_2}}\max_{t\in\R}I_{\beta_2}(t\s(u,v))=m_{\beta_2}(a_1,a_2),$$
so $m_\beta(a_1,a_2)$ is decreasing with respect to $\beta\ge0$.

\vskip0.08in
\item[(2)] Let $l=\min\{m(a_1,0),m(0,a_2)\}$. We first prove $m_0(a_1,a_2)\ge l$. Suppose $0<m_0(a_1,a_2)<l$. Then by Lemma \ref{201906-6cpt}
and Lemma \ref{201909-6PSseq}, we can find a sequence $(u_n,v_n)\to(u_0,v_0)$ in $H$,  where $(u_0,v_0)$ attains the minimum problem
 $m_0(a_1,a_2)$. Since $\beta=0$, the system \eqref{201905-3} is given by two uncoupled equations and both $u_0$ and $v_0$ are
positive radial solutions. By Lemma \ref{201905-nlses}, we have
	\[l>m_0(a_1,a_2)=I_0(u_0,v_0)=m(a_1,0)+m(0,a_2)>l,\]
a contradiction.

	Now we prove $m_0(a_1,a_2)\le l$. Let $u$ be the unique positive solution of \eqref{201905-nlse}
with parameters $p,\mu_1,a_1$ and $v$ be the unique positive solution of \eqref{201905-nlse} with parameters $q,\mu_2,a_2$.
So $(u,v)\in S_{a_1}\times S_{a_2}$ and $(u,s\s v)\in S_{a_1}\times S_{a_2}$ for any $s\in\R$. Let $t_s=t_{(u,s\s v)}$, there is 
\bg
	0&=P_0(t_s\s(u,s\s v))\\
	 &=e^{2t_s}\int_{\RN}|\nabla u|^2+e^{2t_s+2s}\int_{\RN}|\nabla v|^2-e^{p\ga_pt_s}\int_{\RN}\mu_1|u|^p-e^{q\ga_q(t_s+s)}\int_{\RN}\mu_2|v|^q,
\end{align*}
which means that
	\[\int_{\RN}|\nabla u|^2+e^{2s}\int_{\RN}|\nabla v|^2\ge e^{(p\ga_p-2)t_s}\int_{\RN}\mu_1|u|^p.\]
Therefore,  $e^{t_s}$ is bounded as $s\to-\iy$. Hence, for any $s\in\R$,
\bg
	m_0(a_1,a_2)&\le I_0(t_s\s(u,s\s v))=E_{p,\mu_1}(t_s\s u)+E_{q,\mu_2}((t_s+s)\s v)\\
		&\le m(a_1,0)+\f{e^{2(t_s+s)}}{2}\int_{\RN}|\nabla u|^2-\f{e^{q\ga_q(t_s+s)}}{q}\int_{\RN}\mu_2 |v|^q.
\end{align*}
Let $s\to-\iy$, we obtain $m_0(a_1,a_2)\le m(a_1,0)$. Similarly we can prove that $m_0(a_1,a_2)\le m(0,a_2)$.
\end{itemize}
\ep

\bl\lab{201909-6upbdd} 	
	For any $a_1,a_2>0$,  we have that 
	\begin{itemize}[fullwidth,itemindent=0em]
	\item[(1)]	there exists a $\beta_0>0$ such that $m(a_1,a_2)<\min\{m(a_1,0),m(0,a_2)\}$ for any $\beta>\beta_0$;
	\item[(2)]	futher, if $r_1,r_2<2$, then $m(a_1,a_2)<\min\{m(a_1,0),m(0,a_2)\}$ for any $\beta>0$.
	\end{itemize}
\el
\bp
\begin{itemize}[fullwidth,itemindent=0em]
\item[(1)] 	Let $u$ be the unique positive solution of \eqref{201905-nlse} with parameters $p,\mu_1,a_1$ and $v$ be the unique positive
solution of \eqref{201905-nlse} with parameters $q,\mu_2,a_2$. It is easy to see that
	\[E_{p,\mu_1}(s\s u)\to0\quad \text{and}\quad E_{q,\mu_2}(s\s v)\to0\quad as\ s\to-\iy.\]
So there exists a $s_0<-1$ which is independent of $\beta$ such that
\be\lab{201909-6.101}
	\max_{s<s_0}I(s\s(u,v))<\max_{s<s_0}E_{p,\mu_1}(s\s u)+E_{q,\mu_2}(s\s v)<\min\big\{m(a_1,0),m(0,a_2)\big\}.
\ee
If $s\ge s_0$, then the intersection term can be bounded from below:
	\[\int_{\RN}|s\s u|^{r_1}|s\s v|^{r_2}=e^{r\ga_rs}\int_{\RN}|u|^{r_1}|v|^{r_2}\ge Ce^{r\ga_rs_0}.\]
As a consequence, we have
\begin{align*}
	\max_{s\ge s_0}I(s\s(u,v))&\le\max_{s\ge s_0}E_{p,\mu_1}(s\s u)+E_{q,\mu_2}(s\s v)-Ce^{r\ga_rs_0}\beta\\
		&\le m(a_1,0)+m(0,a_2)-Ce^{r\ga_rs_0}\beta,
\end{align*}
and the last term is strictly smaller than $\min\big\{m(a_1,0),m(0,a_2)\big\}$ provided $\beta$ is sufficiently large.

\vskip0.08in
\item[(2)]	Let $(u,\la_0)\in S_{a_1}\times \R^+$ be the unique positive solution in Lemma \ref{201905-nlses} with parameters $p,\mu_1,a_1$. Since $r_2<2$, we take
a $m\in(\f{N}{2}-\f{2}{r_2},\f{N}{2}-1)$ and $v(x)=c\f{\vp(x)}{|x|^m}$ with
	$$\vp(x)\in C_0^\iy(B_2(0)),\quad 0\le\vp(x)\le1,\quad\vp(x)=1\ \text{in}\ B_1(0).$$
Then $v\in H$ and we choose a suitable $c$ such that $v\in S_{a_2}$. Therefore $(u,s\s v)\in S_{a_1}\times S_{a_2}$ for
any $s\in\R$. Let
	\[\al(s)=\int_{\RN}|u|^{r_1}|s\s v|^{r_2}=C_0e^{(\f{N}{2}-m)r_2s}\int_{\RN}u^{r_1}(x)\f{\vp^{r_2}(e^sx)}{|x|^{mr_2}}.\]
As in Lemma \ref{201906-4upbdd}, we can prove that
	\[\al(s)=C_0e^{\theta s}(C_1+o(1))>\f{C_0C_1}{2}e^{\theta s} \quad \text{as}\ s\to-\iy.\]
where $C_1=\int_{\RN}\f{u^{r_1}(x)}{|x|^{mr_2}}\in(0,\iy)$ and $\theta=(\f{N}{2}-m)r_2\in(1,2)$.
Now let $t_s=t_{(u,s\s v)}$, then
\be\lab{201909-6.102}
	\begin{aligned}
	0&=P_0(t_s\s(u,s\s v))\\
	 &=e^{2t_s}\int_{\RN}|\nabla u|^2+e^{2t_s+2s}\int_{\RN}|\nabla v|^2-e^{p\ga_pt_s}\int_{\RN}\mu_1|u|^p\\
	 &\quad\ -e^{q\ga_q(t_s+s)}\int_{\RN}\mu_2|v|^q-\beta r\ga_re^{r\ga_rt_s}\al(s),
	\end{aligned}
\ee
from which we obtain that there exists $C_2,C_3>0$ such that
	\[C_2\le e^{t_s}\le C_3\quad \text{as}\ s\to-\iy.\]
Without loss of generality, we assume $e^{t_s}\to l>0$ as $s\to-\iy$, then let $s\to-\iy$ in \eqref{201909-6.102}, we conclude
	$$l^2\int_{\RN}|\nabla u|^2-l^{p\ga_p}\int_{\RN}\mu_1 |u|^p=0,$$
which menas $l=1$. Therefore
\begin{align*}
	m(a_1,a_2)&\le I(t_s\s(u,s\s v))\\
		&=E_{p,\mu_1}(t_s\s u)+\f{e^{2(t_s+s)}}{2}\int_{\RN}|\nabla v|^2-\f{e^{q\ga_q(t_s+s)}}{q}\int_{\RN}\mu_2|v|^q-\beta e^{r\ga_r t_s}\al(s)\\
		&< m(a_1,0)+\f{e^{2(t_s+s)}}{2}\int_{\RN}|\nabla v|^2-\f{e^{q\ga_q(t_s+s)}}{q}\int_{\RN}\mu_2|v|^q-\beta e^{r\ga_r t_s}\f{C_0C_1}{2}e^{\theta s},\\
\end{align*}
from which, we see for sufficiently small $s\ll-1$, there holds $m(a_1,a_2)<m(a_1,0)$. Similarly we can prove $m(a_1,a_2)<m(0,a_2)$.
\end{itemize}
\ep

\bp[Proof of the Theorem \ref{thm4}]
	The proof is finished when we combine Lemma \ref{201906-6cpt}, Lemma \ref{201909-6PSseq} and Lemma \ref{201909-6upbdd}.
\ep

\vskip0.26in
\appendix
\renewcommand{\appendixname}{Appendix \Alph{section}}
\section{Proof of Lemma \ref{201905-4h} and Remark \ref{tem1}}

\bp[Proof of Lemma \ref{201905-4h}]
Since the monotonicity of $h(t)$ will be strongly affected by the comparision of $p,q$ and $r$, we need to divide the proof into 
four different situations.
\begin{itemize}[fullwidth,itemindent=0em]
\item[{\it Case 1: $p\le r<\bar p$.}]
We have $p\ga_p\le r\ga_r<2<q\ga_q$ and
	$$h'(t)=t^{p\ga_p-1}(t^{2-p\ga_p}-\D_1\beta r\ga_rt^{r\ga_r-p\ga_p}-\D_2p\ga_p-\D_3q\ga_qt^{q\ga_q-p\ga_p}).$$
Denote $g(t)=t^{2-p\ga_p}-\D_1\beta r\ga_rt^{r\ga_r-p\ga_p}-\D_3q\ga_qt^{q\ga_q-p\ga_p}$, there are
	$$h'(t)=t^{p\ga_p-1}\big(g(t)-\D_2p\ga_p\big),$$
	$$g'(t)=t^{r\ga_r-p\ga_p-1}\big[(2-p\ga_p)t^{2-r\ga_r}-\D_1\beta r\ga_r(r\ga_r-p\ga_p)-\D_3q\ga_q(q\ga_q-p\ga_p)t^{q\ga_q-r\ga_r}\big].$$
Let $f(t)=(2-p\ga_p)t^{2-r\ga_r}-\D_3q\ga_q(q\ga_q-p\ga_p)t^{q\ga_p-r\ga_r}$, then
	$$g'(t)=t^{r\ga_r-p\ga_p-1}\big[f(t)-\D_1\beta r\ga_r(r\ga_r-p\ga_p)\big],$$
	$$f'(t)=t^{1-r\ga_r}\big[(2-p\ga_p)(2-r\ga_r)-\D_3q\ga_q(q\ga_q-p\ga_p)(q\ga_q-r\ga_r)t^{q\ga_q-2}\big].$$
Since $p\ga_p\le r\ga_r<2<q\ga_q$, we get $f(+\iy)=g(+\iy)=h(+\iy)=-\iy$,$f(0+)=0^+,h(0+)=0^-$ and
	$$g(0+)=\begin{cases}0^-,\quad p<r,\\ -\D_1\beta r\ga_r<0,\quad p=r. \end{cases}$$
We see that $f(t)$
has a unique critical point $\bar t$ in $(0,+\iy)$ satisfying
\be\lab{201905-4t}
\bar t^{q\ga_q-2}=\f{2-p\ga_p}{q\ga_q-p\ga_p}\f{2-r\ga_r}{q\ga_q-r\ga_r}\f{1}{\D_3q\ga_q}.
\ee
Moreover, if
\be\lab{201905-4hgf}
f(\bar t)>\D_1\beta r\ga_r(r\ga_r-p\ga_p)),\quad g(\bar t)>\D_2p\ga_p,\quad h(\bar t)>0,
\ee
then the function $h(t)$ has exactly two critical points, one is a local minimum at negitive level, the other one is a global
maximum at positive level. Futher, there exists $0<R_0<R_1$ such that $h(R_0)=h(R_1)=0$, and $h(t)>0$ iff $t\in(R_0,R_1)$.
On the other hand, from the definitions of $f(t),g(t)$ and $h(t)$, we can check that \eqref{201905-4hgf} is equivalent to
\be\lab{201905-4hgfe}
\bcs
(2-p\ga_p)\bar t^2>\D_1\beta r\ga_r(r\ga_r-p\ga_p))\bar t^{r\ga_r}+\D_3q\ga_q(q\ga_q-p\ga_p)\bar t^{q\ga_q},\\
\bar t^2>\D_1\beta r\ga_r\bar t^{r\ga_r}+\D_2p\ga_p\bar t^{p\ga_p}+\D_3q\ga_q\bar t^{q\ga_q},\\
\f{1}{2}\bar t^2>\D_1\beta \bar t^{r\ga_r}+\D_2\bar t^{p\ga_p}+\D_3\bar t^{q\ga_q}.
\ecs\ee
Substituting \eqref{201905-4t} into \eqref{201905-4hgfe}, we obtain a constant $C>0$ such that if
	$$\D_1\beta\D_3^{\f{2-r\ga_r}{q\ga_q-2}}+\D_2\D_3^{\f{2-p\ga_p}{q\ga_q-2}}<C,$$
then \eqref{201905-4hgfe} holds, which menas \eqref{201905-4hgf} holds. It follows from the definitions of $\D_1,\D_2$ and $\D_3$ that we can immediately obtain a constant $\al_1$
with the required properties.

\vskip0.08in
\item[{\textit Case 2: $r<p<\bar p$.}]
	If we exchange the roles played by $\D_2t^{p\ga_p}$ and $\D_1t^{r\ga_r}$, then we can get the constant $\al_1$ as {\it Case 1}.

\vskip0.08in
\item[{\textit Case 3: $r=\bar p$.}]
We first suppose $\al_1<\frac{1}{4}$, then $\delta:=\f{1}{2}-\D_1\beta\in(\f{1}{4},\f{1}{2})$ when $\D_1\beta<\al_1$ and $h(t)$ turns to be
	$$h(t)=\delta t^2-\D_2t^{p\ga_p}-\D_3t^{q\ga_q}.$$
Taking a similar argument as in {\it Case 1}, we can prove the existence of the constant $\al_1$.

\vskip0.08in
\item[{\textit Case 4: $r>\bar p$.}]
Note that in this case $p\ga_p<2<r\ga_r,q\ga_q$. Similarly we have
	$$h'(t)=t^{p\ga_p-1}(t^{2-p\ga_p}-\D_1\beta r\ga_rt^{r\ga_r-p\ga_p}-\D_2p\ga_p-\D_3q\ga_qt^{q\ga_q-p\ga_p}).$$
Denote $g(t)=t^{2-p\ga_p}-\D_1\beta r\ga_rt^{r\ga_r-p\ga_p}-\D_3q\ga_qt^{q\ga_q-p\ga_p}$, there are
	$$h'(t)=t^{p\ga_p-1}(g(t)-\D_2p\ga_p),$$
	$$g'(t)=t^{1-p\ga_p}\big[2-p\ga_p-\D_1\beta r\ga_r(r\ga_r-p\ga_p)t^{r\ga_r-2}-\D_3q\ga_q(q\ga_q-p\ga_p)t^{q\ga_q-2}\big].$$
We see that $g(t)$ has a unique critical point $\bar t$ in $(0,+\iy)$ and
\be\lab{201905-44t}
(2-p\ga_p)\bar t^2=\D_1\beta r\ga_r(r\ga_r-p\ga_p)\bar t^{r\ga_r}+\D_3q\ga_q(q\ga_q-p\ga_p)\bar t^{q\ga_p}.
\ee
In particular, if
\be\lab{201905-44hg}
g(\bar t)>\D_2p\ga_p,\quad h(\bar t)>0,
\ee
then $h(t)$ has exactly two critical points: one is a local minimum at a negitive level, the other on is a global maximum at positive
level. Futher, there exist $0<R_0<R_1$ suct that $h(R_0)=h(R_1)=0$, and $h(t)>0$ iff $t\in(R_0,R_1)$. On the other hand, \eqref{201905-44hg} is
equivalent to
\be\lab{201905-44hge}
\bcs
\bar t^2>\D_1\beta r\ga_r\bar t^{r\ga_r}+\D_2p\ga_p\bar t^{p\ga_p}+\D_3q\ga_q\bar t^{q\ga_q},\\
\f{1}{2}\bar t^2>\D_1\beta \bar t^{r\ga_r}+\D_2\bar t^{p\ga_p}+\D_3\bar t^{q\ga_q}.
\ecs\ee
We observe that if
	$$\bar t>\bar s:=\big(2\D_2\min\big\{\f{r\ga_r-2}{r\ga_r-p\ga_p},\f{q\ga_q-2}{q\ga_q-p\ga_p}\big\}\big)^{\f{1}{2-p\ga_p}},$$
then there are
	$$\begin{aligned}
	&\quad\D_1\beta r\ga_r\bar t^{r\ga_r}+\D_2p\ga_p\bar t^{p\ga_p}+\D_3q\ga_q\bar t^{q\ga_q}\\
	&\le\max\big\{\frac{1}{r\ga_r-p\ga_p},\f{1}{q\ga_q-p\ga_p}\big\}(2-p\ga_p)\bar t^2+\D_2q\ga_q\bar s^{p\ga_p-2}\bar t^2\\
	&<\bar t^2,
	\end{aligned}$$
and 
	$$\D_1\beta \bar t^{r\ga_r}+\D_2\bar t^{p\ga_p}+\D_3\bar t^{q\ga_q}<\f{1}{2}\bar t^2.$$
So we just need to guarantee $\bar t>\bar s$. Note that there exists a constant $C>0$ such that
	$$(2-p\ga_p)\bar s^2>\D_1\beta r\ga_r(r\ga_r-p\ga_p)\bar s^{r\ga_r}+\D_3q\ga_q(q\ga_q-p\ga_p)\bar s^{q\ga_p}$$
as long as
	$$\D_1\beta\D_2^{\f{r\ga_r-2}{2}-p\ga_p}+\D_3\D_2^{\f{q\ga_q-2}{q\ga_q-p\ga_p}}<C,$$
then $\bar t>\bar s$ because of $q\ga_q,r\ga_q>2$. Finally, completely analogous to {\it Case 1}, we get the constant $\al_1$ with the
required properties.
\end{itemize}
\ep

\bp[Proof of Remark \ref{tem1}] The proof is similar to the one of Lemma \ref{201905-4h}, and we disscuss in four different cases again.
\begin{itemize}[fullwidth,itemindent=0em]
	\item[{\it Case 1: $p\le r<\bar p$.}]	In this case, we have $p\ga_p\le r\ga_r<2<q\ga_q$ and
	$$l'(t)=t^{p\ga_p-1}\sbr{2at^{2-p\ga_p}-br\ga_rt^{r\ga_r-p\ga_p}-cp\ga_p-dq\ga_qt^{q\ga_q-p\ga_p}}.$$
	We denote $g(t)=2at^{2-p\ga_p}-br\ga_rt^{r\ga_r-p\ga_p}-dq\ga_qt^{q\ga_q-p\ga_p}$, and hence
		$$l'(t)=t^{p\ga_p-1}\sbr{g(t)-cp\ga_p},$$
		$$g'(t)=t^{r\ga_r-p\ga_p-1}\mbr{(2-p\ga_p)2at^{2-r\ga_r}-br\ga_r(r\ga_r-p\ga_p)-dq\ga_q(q\ga_q-p\ga_p)t^{q\ga_q-r\ga_r}}.$$
	Now let $f(t)=(2-p\ga_p)2at^{2-r\ga_r}-dq\ga_q(q\ga_q-p\ga_p)t^{q\ga_q-r\ga_r}$, then
		$$g'(t)=t^{r\ga_r-p\ga_p-1}\mbr{f(t)-br\ga_r(r\ga_r-p\ga_p)},$$
		$$f'(t)=t^{1-r\ga_r}\mbr{(2-p\ga_p)(2-r\ga_r)2a-dq\ga_q(q\ga_q-p\ga_p)(q\ga_q-r\ga_r)t^{q\ga_q-2}}.$$
	We see that $f(t)$ has only one critical point $\bar t$ in $(0,+\iy)$, which is also a maximum point, and that $f(t)$ is strictly increasing 
	in $(0,\bar t)$ and is strictly decreasing in $(\bar t,+\iy)$. To obtain the monotonicity of $g(t)$, we need to compare the value of $f(\bar t)$ 
	and $br\ga_r(r\ga_r-p\ga_p)$. If $br\ga_r(r\ga_r-p\ga_p)\ge f(\bar t)=\max_{t>0}f(t)$, then $g'(t)\le0$ and $g(t)$ is strictly decreasing in $(0,+\iy)$. 
	Since 
	\be\lab{appendix1}
		g(0+)=\left\{\begin{aligned} &0^-, &p<r,\\ &-br\ga_r, &p=r, \end{aligned}\right.
	\ee
	we have that $g(t)<0<cp\ga_p$, and hence $l'(t)<0$, which means that $l(t)$ has no critical points in $(0,+\iy)$. If $br\ga_r(r\ga_r-p\ga_p)<f(\bar t)=\max_{t>0}f(t)$, 
	then by $f(0)=0$, $f(+\iy)=-\iy$, there exist two constants $0\le t_1<\bar t<t_2$ such that $f(t_1)=f(t_2)=br\ga_r(r\ga_r-p\ga_p)$. 
	So $g(t)$ is strictly decreasing for $0<t<t_1$ and $t>t_2$, and is strictly increasing for $t_1<t<t_2$. It follows form \eqref{appendix1} 
	that $g(t)=cp\ga_p$ has at most two solutions in $(0,+\iy)$, which implies that $l(t)$ has at most two critical points in $(0,+\iy)$.

	\vskip 0.08in
	\item[{\it Case 2: $r\le p<\bar p$.}]	As in Lemma \ref{201905-4h}, if we exchange the roles of $bt^{r\ga_r}$ and $ct^{p\ga_p}$, 
	then we can conclude that $l(t)$ has at most two critical points in $(0,+\iy)$.

	\vskip 0.08in
	\item[{\it Case 3: $r=\bar p$.}]	In this case $p\ga_p<r\ga_r=2<q\ga_q$, and hence $l(t)$ becomes 
		$$l(t)=(a-b)t^2-ct^{p\ga_p}-dt^{q\ga_q}.$$
	If $a<b$, we see that $l(t)$ is strictly decreasing and has no critical points. Now suppose $a>b$, then according to {\it Case 1} 
	with $p=r$, we conclude that $l(t)$ has at most two critical points in $(0,+\iy)$.

	\vskip 0.08in
	\item[{\it Case 4: $r>\bar p$.}]	Note that in this case $p\ga_p<2<r\ga_r,q\ga_q$. We have
		$$h'(t)=t^{p\ga_p-1}\sbr{2at^{2-p\ga_p}-br\ga_rt^{r\ga_r-p\ga_p}-cp\ga_p-dq\ga_qt^{q\ga_q-p\ga_p}}.$$
	By denoting $g(t)=2at^{2-p\ga_p}-br\ga_rt^{r\ga_r-p\ga_p}-dq\ga_qt^{q\ga_q-p\ga_p}$, there are
		$$h'(t)=t^{p\ga_p-1}\sbr{g(t)-cp\ga_p},$$
		$$g'(t)=t^{1-p\ga_p}\mbr{(2-p\ga_p)2a-br\ga_r(r\ga_r-p\ga_p)t^{r\ga_r-2}-dq\ga_q(q\ga_q-p\ga_p)t^{q\ga_q-2}}.$$
	Wesee that $g(t)$ has a unique critical point $\bar t$ in $(0,+\iy)$, which is also a maximum point, and $g(t)$ is is strictly increasing 
	in $(0,\bar t)$ and is strictly decreasing in $(\bar t,+\iy)$. To obtain the monotonicity of $l(t)$, we need to compare the value of $g(\bar t)$ 
	and $cp\ga_p$. If $cp\ga_p\ge g(\bar t)=\max_{t>0}g(t)$, then $l'(t)<0$ and $l(t)$ has no critical points in $(0,+\iy)$. If $cp\ga_p<g(\bar t)=\max_{t>0}g(t)$, 
	then \eqref{appendix1} implies that $l'(t)=0$ has atmost two soultions in $(0,+\iy)$, that is $l(t)$ has at most two critical points in $(0,+\iy)$.
\end{itemize}
\ep

\section{A regularity result}
We give a proof of the following facts, which  is  probably  known, but for which we can not find a reference.
\bl\lab{201905-regularity}
	Suppose $\om$ is a domain in $\RN(N\ge 3)$ and $(u,v)\in H_0^1(\om)\times H_0^1(\om)$ is a nonnegative weak solution of
	$$
	\bcs
	-\dl u=f(x,u,v),\\
	-\dl v=g(x,u,v),
	\ecs\quad in\ \om
	$$
	where $f(x,u,v),g(x,u,v):\om\times \R^2\ra \R$ are Carath\'eodory functions satisfying
		$$|f(x,u,v)|+|g(x,u,v)|\le C(|u|+|v|+|u|^{2^*-1}+|v|^{2^*-1}),$$
	for some constant $C>0$. Then $(u,v)$ is a smooth solution.
\el
\bp
We prove that $u,v\in L^p(\om)$ for any $p<\iy$ using Moser iteration, then elliptic regularity theory means that $u,v$ are smooth
functions. Choose $s\ge0$ such that $u,v\in L^{2(s+1)}(\om)$. We shall prove that $u\in L^{2^*(s+1)}(\om)$ so that an obvious
bootstrap argument proves the assertion. Choose $L>0$ and set
	$$\psi=\min\big\{(u+v)^s,L\big\},\ \phi=(u+v)\psi^2,\ \om_L=\big\{x\in\RN:(u(x)+v(x))^s\le L\big\}.$$
In what follows we denote by $C$ various constants independent on $L$. We have
	$$\nabla[(u+v)\psi]=(1+s\chi_{\om_L})\psi\nabla (u+v),$$
	$$\nabla\phi=(1+2s\chi_{\om_L})\psi^2\nabla (u+v),$$
and $\phi\in H_0^1(\om)$. Therefore, we obtain
\begin{align*}
\int_{\om}|\nabla(u+v)|^2\psi^2&\le C\int_{\om}\nabla(u+v)\cdot\nabla\phi=C\int_{\om}[f(x,u,v)+g(x,u,v)]\phi\\
	&\le C\int_{\om}(|u|+|v|+|u|^{2^*-1}+|v|^{2^*-1})\phi\\
	&\le C\int_{\om}(|u|+|v|)^{2(s+1)}+(|u|+|v|)^{2^*-2}[(|u|+|v|)\psi]^2\\
	&\le C(1+\int_{\om}w[(|u|+|v|)\psi]^2),
\end{align*}
where $w(x)=(|u|+|v|)^{2^*-2}\in L^{\f{N}{2}}(\om)$. Then we obtain
\begin{align*}
\int_{\om}|\nabla[(u+v)\psi]|^2&\le C\int_{\om}|\nabla(u+v)|^2\psi^2\le C(1+\int_{\om}w[(|u|+|v|)\psi]^2)\\
	&\le C(1+K\int_{|w|\le K}(|u|+|v|)^{2(s+1)}+\int_{|w|>K}w[(|u|+|v|)\psi]^2))\\
	&\le C(1+K+(\int_{|w|>K}w\f{N}{2})^{\f{2}{N}}(\int_\om[(u+v)\psi]^{2^*})^{\f{2}{2^*}})\\
	&\le C(1+K)+\e_K\int_{\om}|\nabla[(u+v)\psi]|^2,
\end{align*}
where $\e_K\ra0$ as $K\ra+\iy$. Choosing $K$ such that $\e_K<\frac{1}{2}$ we arrive at
	$$\int_{\om_L}|\nabla(u+v)^{s+1}|^2=\int_{\om_L}|\nabla[(u+v)\psi]|^2\le C.$$
Letting $L\ra+\iy$, we get $u^{s+1},v^{s+1}\in H^1(\om)$, hence $u\in L^{2^*(s+1)}(\om)$.
\ep
\bc\lab{201905-Ac}
	Any nonnegative solution of \eqref{201905-2} is smooth solution.
\ec
\bp
In this case, $\om=\RN$ and
	$$f(x,u,v)=-\la_1u+\mu_1 |u|^{p-2}u+\beta r_1|u|^{r_1-2}|v|^{r_2}u$$
	$$g(x,u,v)=-\la_2v+\mu_2 |v|^{q-2}v+\beta r_2|u|^{r_1}|v|^{r_2-2}v,$$
then by Young inequality we have
	$$\begin{aligned}
	|f(x,u,v)|+|g(x,u,v)|&\le C(|u|+|v|+|u|^{p-1}+|v|^{q-1}+|u|^{r-1}+|v|^{r-1})\\
		&\le C(|u|+|v|+|u|^{2^*-1}+|v|^{2^*-1}).
	\end{aligned}$$
Then from Lemma \ref{201905-regularity}, we obtain any nonnegative solution of \eqref{201905-2} is smooth.
\ep

\vskip 0.2in\noindent
{\bf Acknowledgements}\\
The authors thank Nicola Soave   for valuable comments for  preparing the current manuscript: he pointed out a gap in Lemma \ref{201906-4cpt} and gave some comments for the Theorem \ref{thm4}.



 \end{document}